\documentclass{elsart}
\usepackage{latexsym,amssymb}
\usepackage{flafter}
\usepackage{graphicx}
\usepackage{color}
\date{21-09-2013}
\newtheorem{theorem}{Theorem}[section]

\newtheorem{lemma}{Lemma}[section]

\begin{document}
\runauthor{Hou-Biao Li, Ming-Yan Song, Er-Jie Zhong, Xian-Ming Gu etc.}
\begin{frontmatter}
\title{Numerical gradient schemes for heat equations based on the collocation polynomial and Hermite interpolation\thanksref{a}}
\thanks[a]{\small Supported by the National Natural Science Foundation of China (11101071, 1117105, 11271001, 51175443)
and the Fundamental Research Funds for China Scholarship Council.}

\author{Hou-Biao Li\corauthref{b}},
\corauth[b]{Corresponding author.}\ead{lihoubiao0189@163.com}
\author{Ming-Yan Song},
\author{Er-Jie Zhong},
\author{Xian-Ming Gu}
\address{School of Mathematics Sciences, University of Electronic Science and Technology of China, Chengdu, 611731, P. R. China}

\begin{abstract}
As is well-known, the advantage of the high-order compact difference scheme (H-OCD) is unconditionally stable and
convergent with the order $O(\tau^2+h^4)$ under the maximum norm. In this article, a new numerical gradient scheme based on the collocation
polynomial and Hermite interpolation is presented. Moreover, the convergence order of this kind of method is also $O(\tau^2+h^4)$ under the discrete maximum norm when the space step size is just twice the one of H-OCD method, which accelerates the computational process and makes the result much smoother to some extent. In addition, some corresponding analyses are made and the Richardson extrapolation technique is also considered in time direction. The results of numerical experiments are also consistent with these theoretical analysis.
\end{abstract}


\end{frontmatter}


\section{Introduction}

\quad Recently, a great deal of efforts have been devoted to the
development of numerical approximation of heat equation problems
(see, \cite{Hu,S,Sun,Zhang}). It is well known that the traditional
numerical schemes have low accuracy, and thus need fine discretization
to obtain desired accuracy, which leads to many computational
challenges due to the prohibitive computer memory and CPU time
requirements (see, \cite{Sun}).

\quad For heat equations, the forward Euler methods, backward Euler
methods and Crank-Nicolson methods were presented in Ref.\cite{S}. In
addition, three layer implicit schemes also appeared in Ref.\cite{Zhang}.
The forward Euler method and backward method only have one-order accuracy
in time and two-order accuracy in space. Also, the forward Euler method
is not stable when $c\tau/h^2 >1/2$. The three layer implicit compact
format can reach $O(\tau^2+h^4)$, but the format is complex. The
Crank-Nicolson method has two-order accuracy in time and space, which is not better one compared to the high-order compact
difference scheme (see, \cite{Sun}) with two-order accuracy in time and
four-order accuracy in space. The high-order compact difference format (H-OCD)
has many advantages such as using less grid backplane points, high accuracy,
unconditionally stability with the convergence order $O(\tau^2+h^4)$ under the maximum norm.
Therefore, this scheme plays more and more important role in the
numerical solution of partial differential equations and the computational
fluid mechanics field (see, \cite{Ratnesh,Sun.Zhang,Tian}). But the amount of its calculation will increase rapidly with
the increase of grid points.

\quad This article will give a new numerical gradient scheme based
on the collocation polynomial and Hermite interpolation to overcome
the above problem on the high-order compact difference scheme. First, we obtain the intermediate points of the mesh-grid points
by cubic and bi-cubic Hermite interpolation. And then, based on these intermediate
points, the new explicit scheme on the gradient of the discrete
solutions of heat equations is deduced, which will greatly reduce the amount of calculation in
the same accuracy as the high-order compact difference schemes.

\quad The outline of the article is organized as follows. In Section $2$,
the compact difference scheme is derived for one-dimensional heat equations, and the numerical gradient method is presented and then its convergence is analyzed in detail. In Section $3$, we generalize the previous one-dimensional numerical gradient scheme to the two-dimensional one, some similar results are obtained. In addition, the Richardson extrapolation on
time term is also considered. Finally, some numerical results are reported in Section $4$.

\section{One-dimensional Numerical Gradient Schemes Based on the Local Hermite Interpolation and Collocation Polynomial}
\quad For the convenience of description, let us firstly consider the one-dimensional case and then generalize them to the two-dimensional one.

\subsection{The High-Order Compact Difference Scheme in One-dimensional Case}

\quad Firstly, let us consider the following
one-dimensional heat equation problem
\begin{equation}\label{eq1}
\left\{ {\begin{array}{*{20}l}
  \frac{{\partial u}}{{\partial t}}(x,t)& =&c\frac{{\partial^2 u}}{{\partial x^2}}(x,t), & (x,t)\in(0,1)\times(0,T],\\
  u(x,0)&=&\varphi(x), &x\in[0,1],\\
  u(0,t)&=&g_{1}(t), \ u(1,t)=g_{2}(t),& t\in(0,T],\\
\end{array}} \right.
\end{equation}
where $T$ is a positive number. Denote $\Omega=(0,1)\times(0,T]$. In
addition, the solution $u(x,t)$ is assumed to be sufficiently smooth
and has the required continuous partial derivative.

\quad Next, let us recall the compact difference scheme,
which has been introduced in Ref. \cite{Sun.Zhang}.

\quad Let $\Omega_h=\{x_j|x_j=jh, 0\leq j\leq N\}$ be a
uniform partition of $[0,1]$ with the mesh size $h=1/N$ and
$\Omega_\tau =\{t_k|t_k=k\tau,\ \ 0\leq k \leq M\}$ is a uniform
partition of $[0,T]$ with the time step size $\tau=T/M$.
We denote $\Omega_{h\tau}=\Omega_h\times\Omega_\tau$. Let
$\{u^k_j|0\leq j\leq N, \ \ 0\leq k\leq M\}$ be a mesh function
defined on $\Omega_{h\tau}$. For convenience, some other notations
are introduced below.
\[
   \begin{array}{lll}
         [u]^k_j&=&u(x_j,t_k),\;\;\; \;\;u^k_j \approx  u(x_j,t_k),\;\;u^{k+\frac{1}{2}}_j=\frac{{u^k_j+u^{k+1}_j}}{{2}},\\
          \delta_t u^{k+\frac{1}{2}}_j&=&u^{k+1}_j-u^k_j,\;\; u^k_{j-\frac{1}{2}}=\frac{{u^k_j+u^k_{j-1}}}{{2}},\\
          \delta_x u^k_{j-\frac{1}{2}}&=&u^k_j-u^k_{j-1},\;\; \delta^2_x u^k_j=u^k_{j-1}-2u^k_j+u^k_{j+1}.
   \end{array}
\]
\quad In addition, we sometimes use the index pair $(j,k)$ to
represent the mesh point $(x_j,t_k)$. In order to obtain the high-order compact difference scheme on the equation (\ref{eq1}), let us firstly
recall the following lemma.

\begin{lemma}(\cite{Sun.Zhang}).\label{lem:2.1}
Suppose $g(x)\in \mathbb{C}^6[x_{i-1},x_{i+1}]$, then
\begin{equation}\label{eq2}
  \begin{array}{l}
       \frac{{1}}{{12}}[g''(x_{i-1})+10g''(x_i)+g''(x_{i+1})]-
       \frac{{1}}{{h^2}}[g(x_{i-1})-2g(x_i)+g(x_{i+1})]=\frac{{h^4}}{{240}}g^6(\omega_i),
  \end{array}
\end{equation}
where $ \omega_i \in (x_{i-1},x_{i+1})$.
\end{lemma}

\quad Next, let us consider the equation (\ref{eq1}) at the point $(x_j,t_{k+\frac{1}{2}})$. Since
\begin{equation}\label{eq3}
      \frac{{\partial u}}{{\partial t}}(x_j,t_{k+\frac{1}{2}})=c\frac{{\partial^2 u}}{{\partial x^2}}(x_j,t_{k+\frac{1}{2}}), \quad 0\leq i\leq N, \quad 0\leq k \leq M-1,
\end{equation}
then, for $g=[g_0, \ g_1, \ ..., \ g_N]$, we introduce the
 operator $\beta$ with the help of Lemma 2.1.
We denote
\begin{equation}
    \beta g_j=\frac{{1}}{{12}}[g_{j-1}+10g_j+g_{j+1}], \quad 1 \leq j \leq N-1.
\end{equation}
By the famous Taylor formula, we have that
\begin{equation}
    \frac{{1}}{{12\tau}}[\delta_tu^{k+\frac{1}{2}}_{j-1}+10\delta_tu^{k+\frac{1}{2}}_j+\delta_tu^{k+\frac{1}{2}}_{j+1}]=\frac{{c}}{{h^2}}\delta^2_x u(x_j,t_{k+\frac{1}{2}})+R^k_j,
\end{equation}
and
\begin{equation}\label{eq10}
    R^k_j=\tau^2\beta r^k_j+\frac{{ch^4}}{{480}}[\frac{{\partial^6 u}}{{\partial x^6}}(\theta^k_j,t_k)+\frac{{\partial^6 u}}{{\partial x^6}}(\theta ^{k+1}_j,t_{k+1})],
\end{equation}
where $r^k_j=\frac{{1}}{{24}}\frac{{\partial^3u}}{{\partial t^3}}(x_j,\xi^k_j)-
 \frac{{c}}{{8}}\frac{{\partial^4u}}{{\partial x^2 \partial t^2}}(x_j,\eta^k_j)$,
 and $\xi^k_j,\eta^k_j,\theta ^k_j,\theta ^{k+1}_j\in (x_{j-1},x_{j+1})$,$1\leq j \leq N-1, 0\leq k \leq M-1$.

\quad Noting the initial and boundary conditions in the equation (\ref{eq1}), we obtain the following high-order compact difference scheme.
\begin{equation}\label{eq7}
       \begin{array}{l}
               (u^{k+1}_{j-1}-u^k_{j-1})+10(u^{k+1}_j-u^k_j)+(u^{k+1}_{j+1}-u^k_{j+1})\\
               =\frac{{6c\tau}}{{h^2}}\delta^2_x(u^k_j+u^{k+1}_j),
      \end{array}
\end{equation}
where $1\leq j\leq N-1, 0\leq k\leq M-1,$
\begin{equation}\label{eq8}
    u^0_j=\varphi(x_j),\quad 0\leq j\leq N,
\end{equation}
\begin{equation}\label{eq9}
    u^k_0=g_1(k\tau), \ \ u^k_N=g_2(k\tau),\quad 0\leq k\leq M.
\end{equation}

Denoting $u^k_h=[u^k_1, \ u^k_2,..., \ u^k_{N-1}]^T, \quad k=0,\ 1,\ 2, ...,\ M-1$, the above equations (\ref{eq7})-(\ref{eq9}) can be written as
\begin{equation}
   (T_{1}-\frac{{6c\tau}}{{h^2}}T_2)u^{k+1}_h=(T_{1}+\frac{{6c\tau}}{{h^2}}T_2)u^k_h+F_0,
\end{equation}
where
\[
T_{1}  = \left( {\begin{array}{*{20}c}
   {10} & 1 & 0 &  \ldots  & 0 & 0  \\
   1 & {10} & 1 &  \ldots  & 0 & 0  \\
    \vdots  &  \vdots  &  \vdots  &  \ddots  &  \vdots  &  \vdots   \\
   0 & 0 & 0 &  \ldots  & 1 & {10}  \\
\end{array}} \right){\kern 1pt} {\kern 1pt} {\kern 1pt} {\kern 1pt} {\kern 1pt} \texttt{and} {\kern 1pt} {\kern 1pt} {\kern 1pt} {\kern 1pt} {\kern 1pt} {\kern 1pt} T_2  = \left( {\begin{array}{*{20}c}
   { - 2} & 1 & 0 &  \ldots  & 0 & 0  \\
   1 & { - 2} & 1 &  \ldots  & 0 & 0  \\
    \vdots  &  \vdots  &  \vdots  &  \ddots  &  \vdots  &  \vdots   \\
   0 & 0 & 0 &  \ldots  & 1 & { - 2}  \\
\end{array}} \right).
\]

In addition, if we denote
$$ C=\max\{\frac{{c}}{{240}}\max_{0\leq x\leq1,0\leq t\leq T}|\frac{{\partial^6 u(x,t)}}{{\partial x^6}}|,\\
  \frac{1}{24}\max_{0\leq x\leq1,0\leq t\leq T}|\frac{{\partial^3u(x,t)}}{{\partial t^3}}|+\frac{{c}}{{8}}\max_{0\leq x\leq1,0\leq t\leq T}|\frac{{\partial^4u(x,t)}}{{\partial x^2 \partial t^2}}|\},
$$
then, according to \cite{Sun.Zhang}, we have
\[
 |R^k_j|\leq C(\tau^2+h^4),\quad 1\leq j \leq N-1, \ \ 0\leq k\leq M-1.
\]
That is, the truncation error of the compact difference scheme (\ref{eq7}) is $O(\tau^2+h^4)$.

\subsection{One-dimensional Numerical Gradient Scheme Based on the Local Hermite Interpolation and Collocation Polynomial}

\quad As is stated in previous Section 1 and Section 2.1, the compact difference method have some advantages. However,
the amount of calculation will be increased rapidly with the increase of mesh-grid
points, see those numerical experiments in Section 4. In order to deal with this problem, next we give a new numerical
gradient scheme based on the collocation polynomial and Hermite interpolation.

\quad Let $U_h$ be the vector space of the grid function on $\Omega_{\tau h}$.
The $u_h$ denotes the discrete solution satisfying the formula (\ref{eq7})-(\ref{eq9}). Denote
\begin{equation}\label{eq110}
     P_j=\frac{{\partial u(x_j,t)}}{{\partial x}},P^k_j=\frac{{\partial u(x_j,t_k)}}{{\partial x}}.
\end{equation}

Our strategy is as follows:

\begin{enumerate}
  \item First, get the values of points $u^T_j$ by H-OCD scheme
(\ref{eq7})-(\ref{eq9});
  \item And then, obtain the formula (see Eq. (\ref{eqli})) of $P_j$ with the help of collocation
polynomials; i.e,
$$P_j=\frac{{1}}{{12h}}[8u(x_{j+1},t)-8u(x_{j-1},t)+u(x_{j-2},t)-u(x_{j+2},t)];$$
  \item Finally, get the values (see Eq. (\ref{eq12})) of intermediate points $u^t_{j+\frac{1}{2}}$ based on the Hermite interpolation; i.e.,
$$u^t_{j+\frac{1}{2}}=\frac{{1}}{{2}}(u^t_j+u^t_{j+1})+\frac{{h}}{{8}}(P_j-P_{j+1}).$$
\end{enumerate}

\quad Thus, combining H-OCD scheme with the above improvements, a new explicit
\textbf{numerical gradient scheme} for the gradient terms of the
discrete solutions of heat equations is deduced, which will greatly reduce the
amount of calculation in the same accuracy with the high-order
compact difference format. Next, let us give the concrete analysis.

\subsubsection{The Local Hermite Interpolation and Refinement in one-dimensional case}

\quad For convenience, we just consider Hermite cubic and bi-cubic
interpolation functions $u_H(x,t)$ on the interval
$[x_j,x_{j+1}]\subset\Omega_h$, and its vertexes are as follows:
\[
   z_1(x_j,t),\; z_2(x_{j+1},t)\in\Omega_{\tau h}.
\]
On the segment $z_1-z_2$, let the cubic interpolation function satisfy the condition
\[
    u_H(z_1)=u(z_1), \; u_H(z_2)=u(z_2),
\]
\[
    (u_H)_x(z_1)=u_x(z_1),\; (u_H)_x(z_2)=u(z_2).
\]
Based on Ref.\cite{T.S}, we can get the Hermite interpolation polynomial as follows
\begin{equation}\label{eq21}
      u_H(x_{j+\frac{1}{2}},t)=\frac{{1}}{{2}}[u(z_1)+u(z_2)]+\frac{{h}}{{8}}[u_x(z_1)-u_x(z_2)],
\end{equation}
where $j=1,2,...,N-1$. The interpolation errors are
\begin{equation}\label{eq11}
   \begin{array}{lll}
      u_H(x_{j+\frac{1}{2}},t)-u(x_{j+\frac{1}{2}},t)
        &=&\frac{{1}}{{4!}}u_{xxxx}(\xi_j)(x_{j+\frac{1}{2}}-x_j)^2(x_{j+\frac{1}{2}}-x_{j+1})^2\\
        &=&\frac{{h^4}}{{24\times16}}u_{xxxx}(\xi_j),\quad j=1, 2,...,N-1,
   \end{array}
\end{equation}
where $\xi_j$ lies between $z_1$ and $z_2$ (see, \cite{T.S}).
So, by (\ref{eq110}),  we have the refinement computation format
\begin{equation}\label{eq12}
      u^t_{j+\frac{1}{2}}=\frac{{1}}{{2}}[u(x_j,t)+u(x_{j+1},t)]+\frac{{h}}{{8}}(P_j-P_{j+1}),\quad j=1, \ \ 2 ,...,\ \ N-2.
\end{equation}
From (\ref{eq11}), we know that the above refinement schemes have the four-order accuracy in space direction.

\subsubsection{The Collocation Polynomial in one-dimensional case}

\quad From (\ref{eq12}), we know that we must obtain the
expression of $P_j$ in order to get the specific formula of the
intermediate points. Here, we choose the collocation polynomial
method. For convenience, we firstly consider the sub-domain
\[
   [x_{j-1},x_{j+1}]\subset\Omega, \; j=1,2,...,N-1.
\]
Then, we denote
\[
   \xi=x-x_j, \; x\in\Omega, \; j=1, 2, ..., N-1.
\]
In order to get the approximation polynomial of $u$, we consider the polynomial space
\begin{equation}
    H_4=span\{1,\ \xi,\ \xi^2, \ \xi^3, \ \xi^4\},
\end{equation}
and the approximation polynomial of $u$
\begin{equation}
    H(\xi)=a_0+a_1\xi+a_2\xi^2+a_3\xi^3+a_4\xi^4.
\end{equation}
Let
\[
   \begin{array}{lll}
        H(x_{j-1})&=&u^t_{j-1}, \quad \quad \; H(x_{j+1})=u^t_{j+1},\\
        c\frac{{\partial^2H(x_{j-1})}}{{\partial x^2}}&=&\frac{{\partial u}}{{\partial t}}(x_{j-1},t), \quad \quad c\frac{{\partial^2H(x_j)}}{{\partial x^2}}=\frac{{\partial u}}{{\partial t}}(x_j,t),\\
        c\frac{{\partial^2H(x_{j+1})}}{{\partial x^2}}&=&\frac{{\partial u}}{{\partial t}}(x_{j+1},t),\quad \quad j=2,\ \ 3, \ \ ,...\ \ N-2.
  \end{array}
  \]
Thus, by (\ref{eq1}) and (\ref{eq110}), the expression of $P_j$ can be described as follows.
\begin{equation}\label{eqli}
    \begin{array}{lll}
          P_j&=&\frac{{1}}{{2h}}[u(x_{j+1},t)-u(x_{j-1},t)]+\frac{{h}}{{12c}}[\frac{{\partial u}}{{\partial t}}(x_{j-1},t)-
               \frac{{\partial u}}{{\partial t}}(x_{j+1},t)]\\
            &=&\frac{{1}}{{12h}}[8u(x_{j+1},t)-8u(x_{j-1},t)+u(x_{j-2},t)-u(x_{j+2},t)], \\ \label{eq13}
          P_1&=&\frac{{1}}{{6h}}[-2g_1(t)-3u(x_1,t)+6u(x_2,t)-u(x_3,t)],\\
          P_{N-1}&=&\frac{{1}}{{6h}}[2g_2(t)+3u(x_{N-1},t)-6u(x_{N-2},t)+u(x_{N-3},t)],
    \end{array}
\end{equation}
where $j=2,\ \ 3,...\ \ N-3$.

So
\begin{equation}
    \begin{array}{lll}
    u^t_{j+\frac{1}{2}}&=&\frac{{1}}{{2}}[u(x_j,t)+u(x_{j+1},t)]+\frac{{h}}{{8}}(P_j-P_{j+1})\\
                       &=&\frac{{1}}{{2}}[u(x_j,t)+u(x_{j+1},t)]+\frac{{1}}{{96}}[8u(x_{j+1},t)-8u(x_{j-1},t)+u(x_{j-2},t)-\\
                        &&u(x_{j+2},t)]-\frac{{1}}{{96}}[8u(x_{j+2},t)-8u(x_j,t)+u(x_{j-1},t)-u(x_{j+3},t)]\\
                       &=&\frac{{1}}{{96}}[56u(x_j,t)+56u(x_{j+1},t)+u(x_{j-2},t)-9u(x_{j+2},t)+u(x_{j+3},t)],\\
    u^t_{\frac{3}{2}}&=&\frac{{1}}{{2}}[u(x_1,t)+u(x_2,t)]+\frac{{h}}{{8}}(P_1-P_2)\\
                     &=&\frac{{1}}{{96}}[50u(x_1,t)+60u(x_2,t)-10u(x_3,t)-u(x_0,t)+u(x_4,t)-4g_1(t)],\\
    u^t_{N-\frac{3}{2}}&=&\frac{{1}}{{2}}[u(x_{N-2},t)+u(x_{N-1},t)]+\frac{{h}}{{8}}(P_{N-2}-P_{N-1})\\
                       &=&\frac{{1}}{{96}}[60u(x_{N-2},t)+50u(x_{N-1},t)-10u(x_{N-3},t)+u(x_{N-4},t)-\\
                       &&u(x_N,t)-4g_2(t)],
    \end{array}
\end{equation}
where $j=2,\ \ 3, \ \ ,...\ \ N-3$.

\quad Next, according to our improvement scheme, let us analyze the
convergence order of this kind of numerical gradient scheme.

\begin{theorem}
If $ u(x,t)\in \mathbb{C}^6_x(\Omega) $ and $ u(x,t)\in \mathbb{C}^3_t(\Omega) $, then
we have
\begin{equation}
        |\frac{{\partial u(x_j,t_k)}}{{\partial x}}-P^k_j|\leq O(h^4),
\end{equation}
where $j=2,\ 3, \ 4,\ ...\ N-2 $.
\end{theorem}

\textbf{Proof.} When $T=k\tau$, by (\ref{eq13}), we know that
\[
   \begin{array}{lll}
     P^k_j&=&\frac{{1}}{{12h}}[8u(x_{j+1},t_k)-8u(x_{j-1},t_k)+u(x_{j-2},t_k)-u(x_{j+2},t_k)]\\
        &=&\frac{{4}}{{3}}\frac{{u(x_{j+1},t_k)-u(x_{j-1},t_k)}}{{2h}}-\frac{{1}}{{3}}\frac{{u(x_{j+2},t_k)-u(x_{j-2},t_k)}}{{4h}}.\\
   \end{array}
\]
So, by the H-OCD method and the energy method with the Sobolev embedding theorem in \cite{Liao},
\[|\frac{{\partial u(x_j,t_k)}}{{\partial x}}-\{\frac{{4}}{{3}}\frac{{u(x_{j+1},t_k)-u(x_{j-1},t_k)}}{{2h}}-
  \frac{{1}}{{3}}\frac{{u(x_{j+2},t_k)-u(x_{j-2},t_k)}}{{4h}}\}|=O(h^4),
\]
where $j=2,\ 3, \ 4,\ ...\ N-2 $. Thus the proof is completed. $\Box$

\quad By the above theorem, we know that the accuracy of the
partial derivative of $u$ (i.e., $P_j$) in space direction is $O(h^4)$ when $T=k\tau$. In
fact, due to (\ref{eq12}), it is easy to prove that the accuracy of
the intermediate points is $O(h^4)$, too. The corresponding analysis is as
follows.
\begin{theorem}
If $ u(x,t)\in \mathbb{C}^6_x(\Omega) $ and $ u(x,t)\in \mathbb{C}^3_t(\Omega) $, and $u(x,t)$ is the exact solution
 of the equation (\ref{eq1}), then
 \begin{equation}
        |u(x_{j+\frac{1}{2}},t_k)-u^k_{j+\frac{1}{2}}|\leq O(\tau^2+h^4),
 \end{equation}
 where $j=2,\ 3,\ ... \, N-3$.
\end{theorem}
\textbf{Proof.} First, note that
\[
\begin{array}{lll}
 |u(x_{j+\frac{1}{2}},t_k)-u^k_{j+\frac{1}{2}}|&=& |u(x_{j+\frac{1}{2}},t_k)-u_H(x_{j+\frac{1}{2}},t_k)+u_H(x_{j+\frac{1}{2}},t_k)-u^k_{j+\frac{1}{2}}|\\
 &\leq & |u(x_{j+\frac{1}{2}},t_k)-u_H(x_{j+\frac{1}{2}},t_k)|+|u_H(x_{j+\frac{1}{2}},t_k)-u^k_{j+\frac{1}{2}}|,
  \end{array}
\]
then, by the Taylor expansion at $(x_{j+\frac{1}{2}},T)$ (Here $T=k\tau$), we have
\[
     u(x_j,t_k)=u(x_{j+\frac{1}{2}},t_k)-\frac{{h}}{{2}}u_x(x_{j+\frac{1}{2}},t_k)+\frac{{h^2}}{{8}}u_{xx}(x_{j+\frac{1}{2}},t_k)-
     \frac{{h^3}}{{48}}u_{xxx}(x_{j+\frac{1}{2}},t_k)+O(h^4),
\]
\[
    u(x_{j+1},t_k)=u(x_{j+\frac{1}{2}},t_k)+\frac{{h}}{{2}}u_x(x_{j+\frac{1}{2}},t_k)+\frac{{h^2}}{{8}}u_{xx}(x_{j+\frac{1}{2}},t_k)+
    \frac{{h^3}}{{48}}u_{xxx}(x_{j+\frac{1}{2}},t_k)+O(h^4).
\]
Thus, by (\ref{eq12}), we can obtain
\[
    \begin{array}{lll}
       u(x_{j+\frac{1}{2}},t_k)&=&\frac{{1}}{{2}}[u(x_j,t_k)+u(x_{j+1},t_k)]-\frac{{h^2}}{{8}}u_{xx}(x_{j+\frac{1}{2}},t_k)+O(h^4)\\
                      &=&\frac{{1}}{{2}}[u(x_j,t_k)+u(x_{j+1},t_k)]-\frac{{h^2}}{{8}}\frac{{\partial}}{{\partial x}}\frac{{\partial u}}{{\partial x}}(x_{j+\frac{1}{2}},t_k)+O(h^4)\\
                      &=&\frac{{1}}{{2}}[u(x_j,t_k)+u(x_{j+1},t_k)]-\frac{{h^2}}{{8}}\frac{{\partial}}{{\partial x}}[\frac{{1}}{{h}}\delta_xu(x_{j+\frac{1}{2}},t_k)+O(h^2)]+O(h^4)\\
                      &=&\frac{{1}}{{2}}[u(x_j,t_k)+u(x_{j+1},t_k)]+\frac{{h^2}}{{8}}\frac{{\partial}}{{\partial x}}
                      \frac{{u(x_j,t_k)-u(x_{j+1},t_k)}}{{h}}+O(h^4)\\
                      &=&\frac{{1}}{{2}}[u(x_j,t_k)+u(x_{j+1},t_k)]+\frac{{h}}{{8}}[u_x(x_j,t_k)-u_x(x_{j+1},t_k)]+O(h^4)\\
                       &=&u_H(x_{j+\frac{1}{2}},t_k)+O(h^4),
        \end{array}
\]
where $j=2,\ 3,\ ... \, N-3$.

Note that there is no change in time direction corresponding to H-OCD method. Therefore, by (\ref{eq21}), (\ref{eq12}) and Theorem 2.1, we have
\[
    \begin{array}{lll}
      u_H(x_{j+\frac{1}{2}},T)- u^T_{j+\frac{1}{2}}&=&\frac{{1}}{{2}}[u(z_1)+u(z_2)]+\frac{{h}}{{8}}[u_x(z_1)-u_x(z_2)]-\\
       &&\frac{{1}}{{2}}[u(x_j,T)+u(x_{j+1},T)]-\frac{{h}}{{8}}(P_j-P_{j+1})\\
       &=&\frac{{h}}{{8}}[u_x(z_1)-u_x(z_2)]-\frac{{h}}{{8}}(P_j-P_{j+1})\\
       &=&\frac{{h}}{{8}}[u_x(z_1)-P_j-u_x(z_2)+P_{j+1}]\\
       &\leq&O(\tau^2+h^4)
        \end{array}
\]

So
\[
\begin{array}{lll}
 |u(x_{j+\frac{1}{2}},t_k)-u^k_{j+\frac{1}{2}}| &\leq &
 |u(x_{j+\frac{1}{2}},t_k)-u_H(x_{j+\frac{1}{2}},t_k)|+|u_H(x_{j+\frac{1}{2}},t_k)-u^k_{j+\frac{1}{2}}|\\
  &\leq &O(h^4)+O(\tau^2+h^4)\\
  &\leq &O(\tau^2+h^4).
  \end{array}
\]

Thus the proof is completed. $\Box$

\subsection{The Richardson Extrapolation on the H-OCD Scheme in one-dimensional case}

\quad For the compact difference H-OCD scheme considered in Section 2.1, the numerical solution and its difference quotient in space
direction are unconditionally convergent with the convergence order $O(\tau^2+h^4)$ under the maximum norm. And the convergence of the difference quotient in space direction may be proved by the energy method with the Sobolev embedding theorem, that is,
\[|\frac{{\partial u(x_j,t_k)}}{{\partial x}}-\{\frac{{4}}{{3}}\frac{{u(x_{j+1},t_k)-u(x_{j-1},t_k)}}{{2h}}-
  \frac{{1}}{{3}}\frac{{u(x_{j+2},t_k)-u(x_{j-2},t_k)}}{{4h}}\}|=O(\tau^2+h^4),
\]
and
\[|u(x_j,t_k)-u^k_j|=O(\tau^2+h^4), \quad 1\leq j\leq N-1,\ \ 1\leq k\leq M.
\]

\quad Next, we consider the Richardson extrapolation
on this H-OCD scheme (\ref{eq7})-(\ref{eq9}) in time direction in order to reduce the total computing time by \cite{H.L.}.

\begin{lemma}(\cite{Z}).\label{lem:4.1}
 Let $\{V^k_j|0\leq j\leq N,0\leq k\leq M\}$ be the solution of the equation below
\begin{equation}\
  \begin{array}{l}
      \frac{{1}}{{\tau}}\delta_tV^{k+\frac{1}{2}}_j-\frac{{a}}{{h^2}}\delta^2_xV^{k+\frac{1}{2}}_j=g^{k+\frac{1}{2}}_j,\quad \quad 0\leq j\leq N-1,0\leq k\leq M-1,\\
      V^0_j=\varphi_j,\quad \quad 0\leq j\leq N-1,\\
      V^k_0=0, V^k_N=0,\quad \quad 0\leq k\leq M.
  \end{array}
\end{equation}
then
\[
  |V^k|_\infty \leq \frac{{1}}{{2}}(|V^0|^2_1+\frac{{\tau}}{{2a}}\sum^{k-1}_{l=0}|g^{l+\frac{1}{2}}|^2)^{\frac{1}{2}},\quad \quad 0\leq k\leq N,
\]
where
\[
|g^{l+\frac{1}{2}}|^2=h\sum^{N-1}_{j=1}(g^{l+\frac{1}{2}}_j)^2.
\]
\end{lemma}

\begin{theorem}
Let $u^k_j(h,\tau)$ be the solution of H-OCD scheme (\ref{eq7})-(\ref{eq9}) with the time step $\tau$ and the space step $h$.
Then
\[|u(x_j,t_k)-[\frac{{4}}{{3}}u^{2k}_j(h,\frac{\tau}{2})-\frac{{1}}{{3}}u^k_j(h,\tau)]|
               =O(\tau^4+h^4),\; 1\leq j\leq N-1,\;1\leq k\leq M.
\]
\end{theorem}
\textbf{Proof.} Let us consider the following initial-boundary problem
\[
     \begin{array}{lll}
  u_t-\triangle u&=&F_u(x,t),\quad (x,t)\in (0,1)\times(0,T],\\
u(0,t)&=&u(1,t)=0,\quad 0\leq t\leq T,\\
u(x,0)&=&0, \quad x\in (0,1)
      \end{array}
\]
with the smooth solution $p(x,t)$, where
\[F_p(x,t)=\frac{{1}}{{24}}\frac{{\partial^3u(x,t)}}{{\partial t^3}}-
      \frac{{c}}{{8}}\frac{{\partial^4u(x,t)}}{{\partial x^2 \partial t^2}}.
\]

By (\ref{eq10}), we know
\[R^k_j=F_p(x_j,t_{k+\frac{1}{2}})\tau^2+O(\tau^4+h^4),\quad 1\leq j\leq N-1,\ \ 0\leq k\leq M-1.
\]
So
\[ \begin{array}{lll}
  &&\frac{{\delta_te^{k+\frac{1}{2}}_j}}{{\tau}}-\frac{{c}}{{h^2}}\delta^2_xe^{k+\frac{1}{2}}_j
          =F_p(x_j,t_{k+\frac{1}{2}})\tau^2+O(\tau^4+h^4),\; 1\leq j\leq N-1,\; 0\leq k\leq M-1,\\
 &&e^0_j=0, \quad 0\leq j\leq N,\\
 &&e^k_0=0, \quad e^k_m=0,\quad 1\leq k\leq M.
\end{array}
\]
Here $e^k_j=u(x_j,t_k)-u^k_j,\;0\leq j\leq N, \; 0\leq k\leq M.$

\quad In addition, according to H-OCD (\ref{eq7})-(\ref{eq9}), we obtain
\[ \begin{array}{lll}
 &&\frac{{\delta_tp^{k+\frac{1}{2}}_j}}{{\tau}}-\frac{{c}}{{h^2}}\delta^2_xp^{k+\frac{1}{2}}_j
       =F_p(x_j,t_{k+\frac{1}{2}}),\; 1\leq j\leq N-1,\; 0\leq k\leq M-1,\\
 &&p^0_j=0, \quad 0\leq j\leq N,\\
 &&p^k_0=0, \quad p^k_m=0,\quad 1\leq k\leq M.
\end{array}
\]

Then
\[u(x_j,t_k)-p^k_j(h,\tau)=O(\tau^2+h^4),\quad 1\leq j\leq N-1,\ \ 0\leq k\leq M.
\]
Denote
\[ r^k_j=e^k_j+\tau^2p^k_j,\quad 1\leq j\leq N,\ \ 0\leq k\leq M.
\]

Combine the above equations, we get
\[ \begin{array}{lll}
&& \frac{{\delta_tr^{k+\frac{1}{2}}_j}}{{\tau}}-\frac{{c}}{{h^2}}\delta^2_xr^{k+\frac{1}{2}}_j
     =O(\tau^4+h^4),\quad 1\leq j\leq N-1,\ \ 0\leq k\leq M-1,\\
&&r^0_j=0, \quad 0\leq j\leq N,\\
&&r^k_0=0, \quad r^k_N=0,\quad 1\leq k\leq M,
\end{array}
\]
and then by Lemma \ref{lem:4.1}, we have
\[
|r^k|_\infty \leq \frac{{1}}{{2}}(|r^0|^2_1+\frac{{\tau}}{{2a}}\sum^{k-1}_{l=0}|O^{l+\frac{1}{2}}(\tau^4+h^4)|^2)^{\frac{1}{2}},\quad \quad 0\leq k\leq N,
\]
where
\[
|O^{l+\frac{1}{2}}(\tau^4+h^4)|^2=h\sum^{N-1}_{j=1}(g^{l+\frac{1}{2}}_j(\tau^4+h^4))^2.
\]

That is
\[ |r^k|_\infty=O(\tau^4+h^4),\quad 1\leq k\leq M,
\]
i.e.,
\[u^k_j(h,\tau)=u(x_j,t_k)+\tau^2p(x_j,t_k)+O(\tau^4+h^4),\quad 1\leq j\leq N-1,\ \ 0\leq k\leq M,
\]
\[u^{2k}_j(h,\frac{{\tau}}{{2}})=u(x_j,t_k)+(\frac{{\tau}}{{2}})^2p(x_j,t_k)+
           O((\frac{{\tau}}{{2}})^4+h^4),\quad 1\leq j\leq N-1,\ \ 0\leq k\leq M.
\]
Then
 \[|u(x_j,t_k)-[\frac{{4}}{{3}}u^{2k}_j(h,\frac{\tau}{2})-\frac{{1}}{{3}}u^k_j(h,\tau)]|
               =O(\tau^4+h^4),\quad 1\leq j\leq N-1,\ \ 1\leq k\leq M.\]
Thus the conclusion is proved. $\Box$

\textbf{Remark 2.1} With the Richardson extrapolation method above, the truncation errors in time direction for H-OCD scheme is $O(\tau^4+h^4)$ in terms of the maximum norm. Similarly, the extrapolation
  $\frac{{16}}{{15}}u^{4k}_{2j}(\frac{{\tau}}{{4}},\frac{{h}}{{2}})-\frac{{1}}{{15}}u^k_j(\tau,h)$ will can obtain the following result for any $1\leq j\leq N-1, 1\leq k\leq M$.
\[
|u(x_j,t_k)-[\frac{{16}}{{15}}u^{4k}_{2j}(\frac{{\tau}}{{4}},\frac{{h}}{{2}})-\frac{{1}}{{15}}u^k_j(\tau,h)]|=O(\tau^4+h^6).
\]

\section{Two-dimensional Numerical Gradient Scheme Based on the Local Hermite Interpolation and Collocation Polynomial}

\subsection{The High-Order Compact Difference Scheme in Two-dimensional Case}
\quad Next, let us generalize the previous one-dimensional H-OCD scheme to the two-dimensional one. Similar to the previous Section 2, the following
two-dimensional heat equation problem is considered
\begin{equation}\label{eq3.1}
\left\{ {\begin{array}{*{20}l}
  \frac{{\partial u}}{{\partial t}}=\frac{{\partial^2 u}}{{\partial x^2}}+\frac{{\partial^2 u}}{{\partial y^2}}, &(x,y,t)\in \Omega\times(0,T],\\
  u(x,y,0)=\varphi(x,y), & (x,y)\in[a,b]\times[c,d], \\
  u(a,y,t)=g_{1}(y,t), \ u(b,y,t)=g_{2}(y,t),&\ \ (y,t)\in[c,d]\times(0,T],\\
   u(x,c,t)=g_{3}(x,t), \ u(x,d,t)=g_{4}(x,t),&\ \ (x,t)\in[a,b]\times(0,T],
\end{array}} \right.
\end{equation}
where $T$ is a positive number. Denote $\Omega=(a,b)\times(c,d)$. In
addition, the solution $u(x,y,t)$ is assumed to be sufficiently smooth
and has the required continuous partial derivative.

\quad Let ${h_x} = \frac{{b - a}}{{{N_x}}},{h_y} = \frac{{d - c}}{{{N_y}}}$ and ${\Omega _h} = \{ \left( {{x_i},{y_j}} \right)|{x_i} = i{h_x},{y_j} = j{h_y},0 \le i,j \le N\} $. When $\tau = T/M$, define ${{{\Omega }}_{{\tau }}} = \{ {{{t}}_{{k}}}|{{{t}}_{{k}}} = {{k\tau }},0 \le {{k}} \le {{M}}\} $ and ${\Omega _{h\tau }} = {\Omega _h} \times {\Omega _\tau }$. In addition, we denote $\{ u_{{{ij}}}^{{k}}|0 \le {{i}}, {{j}} \le {{N}},0 \le {\rm{k}} \le {\rm{M}}\}$ the mesh function defined on $\Omega_{h\tau}$. Moreover, some other notations
are introduced below.
\[
\begin{array}{lll}
\left[u \right]_{ij}^k = u\left( {{x_i},{y_j},{t_k}} \right),u_{ij}^k \approx u\left( {{x_i},{y_j},{t_k}} \right);\\
\delta _x^2u_{ij}^k = u_{i - 1,j}^k - 2u_{ij}^k + u_{i + 1,j}^k,\;\; \delta _y^2u_{ij}^k = u_{i,j - 1}^k - 2u_{ij}^k + u_{i,j + 1}^k.
\end{array}
\]
For convenience, define the operators ${D_x} = \frac{d}{{dx}},{D_y} = \frac{d}{{dy}}$, ${\rm{E}} = \mathop \sum \limits_{k = 0}^\infty  \frac{1}{{k!}}{\left( {hD} \right)^k} = {e^{hD}}$ and its inverse operator ${E^{-1}} = {e^{ - hD}}$. Obviously
$${\delta ^2} = {E^{ - 1}} - 2 + E = {e^{ - hD}} - 2 + {e^{hD}} = {h^2}{D^2} + 1/12{h^4}{D^4} + O({h^6}).$$
Note that ${\delta ^2} = {h^2}{D^2} + O\left( {{h^4}} \right)$, therefore
\begin{equation}\label{eq3.7}
{\delta ^2} = {h^2}{D^2} + \frac{1}{{12}}{h^2}{D^2}\left[ {{\delta ^2} - O\left( {{h^4}} \right)} \right] + O\left( {{h^6}} \right) = \left( {1 + \frac{1}{{12}}{\delta ^2}} \right){h^2}{D^2} + O\left( {{h^6}} \right).
\end{equation}
i.e.,
$${(1 + \frac{1}{{12}}{\delta ^2})^{ - 1}}{\delta ^2} = {h^2}{D^2} + O\left( {{h^6}} \right).$$
Apply (\ref{eq3.7}) to (\ref{eq3.1}), we obtain that
\begin{equation}\label{eq3.11}
\left( {1 + \frac{1}{{12}}\delta _x^2 + \frac{1}{{12}}\delta _y^2} \right)\frac{{\partial u}}{{\partial t}} = \frac{{\left( {1 + \frac{1}{{12}}\delta _y^2} \right)}}{{h_x^2}}\delta _x^2u + \frac{{\left( {1 + \frac{1}{{12}}\delta _x^2} \right)}}{{h_y^2}}\delta _y^2u + O\left( {{h^4}} \right),
\end{equation}
where $O\left( {{h^4}} \right) = O\left( {h_x^4 + h_y^4} \right)$. In addition, according to the Crank-Nicolson difference scheme \cite{Hu}, we further have that
$$\left( {1 + \frac{1}{{12}}\delta _x^2 + \frac{1}{{12}}\delta _y^2} \right)\frac{{{u^{n + 1}} - {u^n}}}{\tau } = \left[ {\frac{{\left( {1 + \frac{1}{{12}}\delta _y^2} \right)}}{{h_x^2}}\delta _x^2 + \frac{{\left( {1 + \frac{1}{{12}}\delta _x^2} \right)}}{{h_y^2}}\delta _y^2} \right]\frac{{{u^{n + 1}} + {u^n}}}{2},$$
which can be written as
\[
\begin{array}{l}
\left( {1 + \frac{1}{{12}}\delta _x^2 + \frac{1}{{12}}\delta _y^2 - \tau \frac{{\left( {1 + \frac{1}{{12}}\delta _y^2} \right)}}{{2h_x^2}}\delta _x^2 - \tau \frac{{\left( {1 + \frac{1}{{12}}\delta _x^2} \right)}}{{2h_y^2}}\delta _y^2} \right){u^{n + 1}}\\
 = \left( {1 + \frac{1}{{12}}\delta _x^2 + \frac{1}{{12}}\delta _y^2 + \tau \frac{{\left( {1 + \frac{1}{{12}}\delta _y^2} \right)}}{{2h_x^2}}\delta _x^2 + \tau \frac{{\left( {1 + \frac{1}{{12}}\delta _x^2} \right)}}{{2h_y^2}}\delta _y^2} \right){u^n}.
\end{array}
\]
Next, in order to the convenience of description, let $h = {h_x} = {h_y}, N = {N_x} = {N_y}$ and define $r = \tau/{{2{h^2}}}$, then the above equation can be reduced to the following discrete form by the initial and boundary conditions
\begin{equation}\label{eq3.15}
\begin{array}{l}
\left( {\frac{2}{3} + \frac{{10}}{3}r} \right)u_{ij}^{n + 1} - \left( {\frac{{2r}}{3} - \frac{1}{{12}}} \right)\left( {u_{i - 1,j}^{n + 1} + u_{i + 1,j}^{n + 1} + u_{i,j - 1}^{n + 1} + u_{i,j + 1}^{n + 1}} \right) - \frac{r}{6}(u_{i - 1,j - 1}^{n + 1} +u_{i - 1,j + 1}^{n + 1}\\
 + u_{i + 1,j - 1}^{n + 1} + u_{i + 1,j + 1}^{n + 1}) = \left( {\frac{{2r}}{3} + \frac{1}{{12}}} \right)\left( {u_{i - 1,j}^n + u_{i + 1,j}^n + u_{i,j - 1}^n + u_{i,j + 1}^n} \right) + \\
\left( {\frac{2}{3} - \frac{{10}}{3}r} \right)u_{ij}^n + \frac{r}{6}\left( {u_{i - 1,j - 1}^n + u_{i - 1,j + 1}^n + u_{i + 1,j - 1}^n + u_{i + 1,j + 1}^n} \right),
\end{array}
\end{equation}
where $i,j = 1,2, \ldots, N - 1$. $u_{ij}^0 = {u_0}\left( {{x_i},{y_j}} \right),u_{0j}^k = {g_1}\left( {{y_j},{t_k}} \right),u_{Nj}^k = {g_2}\left( {{y_j},{t_k}} \right)$, $u_{i0}^k = {g_3}\left( {{x_i},{t_k}} \right),u_{iN}^k = {g_4}\left( {{x_i},{t_k}} \right),\; i,j = 1,2, \ldots ,N - 1,k = 1,2, \ldots ,M.$
The concrete computation process of the above discrete scheme may be described as follows.
\begin{equation}\label{eq3.17}
\begin{array}{l}
\left( {\begin{array}{*{20}{c}}
{\begin{array}{*{20}{c}}
{\begin{array}{*{20}{c}}
{{A_1}}\\
{ - {B_1}}
\end{array}}&{\begin{array}{*{20}{c}}
{ - {B_1}}\\
{{A_1}}
\end{array}}
\end{array}}&{\begin{array}{*{20}{c}}
 \cdots \\
 \cdots
\end{array}}&{\begin{array}{*{20}{c}}
{\begin{array}{*{20}{c}}
0&0
\end{array}}\\
{\begin{array}{*{20}{c}}
0&0
\end{array}}
\end{array}}\\
{\begin{array}{*{20}{c}}
 \vdots & \vdots
\end{array}}& \ddots &{\begin{array}{*{20}{c}}
 \vdots & \vdots
\end{array}}\\
{\begin{array}{*{20}{c}}
0&0
\end{array}}& \cdots &{\begin{array}{*{20}{c}}
{ - {B_1}}&{{A_1}}
\end{array}}
\end{array}} \right)\left( {\begin{array}{*{20}{c}}
{\begin{array}{*{20}{c}}
{u_{h1}^{n + 1}}\\
{u_{h2}^{n + 1}}
\end{array}}\\
{\begin{array}{*{20}{c}}
 \vdots \\
{u_{hN - 1}^{n + 1}}
\end{array}}
\end{array}} \right) = \left( {\begin{array}{*{20}{c}}
{\begin{array}{*{20}{c}}
{\begin{array}{*{20}{c}}
{{A_2}}\\
{{B_2}}
\end{array}}&{\begin{array}{*{20}{c}}
{{B_2}}\\
{{A_2}}
\end{array}}
\end{array}}&{\begin{array}{*{20}{c}}
 \cdots \\
 \cdots
\end{array}}&{\begin{array}{*{20}{c}}
{\begin{array}{*{20}{c}}
0&0
\end{array}}\\
{\begin{array}{*{20}{c}}
0&0
\end{array}}
\end{array}}\\
{\begin{array}{*{20}{c}}
 \vdots & \vdots
\end{array}}& \ddots &{\begin{array}{*{20}{c}}
 \vdots & \vdots
\end{array}}\\
{\begin{array}{*{20}{c}}
0&0
\end{array}}& \cdots &{\begin{array}{*{20}{c}}
{{B_2}}&{{A_2}}
\end{array}}
\end{array}} \right)\left( {\begin{array}{*{20}{c}}
{\begin{array}{*{20}{c}}
{u_{h1}^n}\\
{u_{h2}^n}
\end{array}}\\
{\begin{array}{*{20}{c}}
 \vdots \\
{u_{hN - 1}^n}
\end{array}}
\end{array}} \right) + \\
 \left( {\begin{array}{*{20}{c}}
{\begin{array}{*{20}{c}}
{U_{h1}^{n + 1}}\\
{U_{h2}^{n + 1}}
\end{array}}\\
{\begin{array}{*{20}{c}}
 \vdots \\
{U_{hN - 1}^{n + 1}}
\end{array}}
\end{array}} \right) + \left( {\begin{array}{*{20}{c}}
{\begin{array}{*{20}{c}}
{{B_1}u_{h0}^{n + 1}}\\
0
\end{array}}\\
{\begin{array}{*{20}{c}}
 \vdots \\
{{B_1}u_{hN}^{n + 1}}
\end{array}}
\end{array}} \right) + \left( {\begin{array}{*{20}{c}}
{\begin{array}{*{20}{c}}
{U_{h1}^n}\\
{U_{h2}^n}
\end{array}}\\
{\begin{array}{*{20}{c}}
 \vdots \\
{U_{hN - 1}^n}
\end{array}}
\end{array}} \right) + \left( {\begin{array}{*{20}{c}}
{\begin{array}{*{20}{c}}
{{B_2}u_{h0}^n}\\
0
\end{array}}\\
{\begin{array}{*{20}{c}}
 \vdots \\
{{B_2}u_{hN}^n}
\end{array}}
\end{array}} \right),
\end{array}
\end{equation}
where
\[\begin{array}{l}
U_{hj}^{n + 1} = {\left[ {\frac{{8r - 1}}{{12}}u_{0j}^{n + 1} + \frac{r}{6}u_{0j + 1}^{n + 1} + \frac{r}{6}u_{0j - 1}^{n + 1},0 \ldots 0,\frac{{8r - 1}}{{12}}u_{Nj}^{n + 1} + \frac{r}{6}u_{N,j + 1}^{n + 1} + \frac{r}{6}u_{N,j - 1}^{n + 1}} \right]^T},\\
U_{hj}^n = {\left[ {\frac{{8r + 1}}{{12}}u_{0j}^n + \frac{r}{6}u_{0j + 1}^n + \frac{r}{6}u_{0j - 1}^n,0 \ldots 0,\frac{{8r + 1}}{{12}}u_{Nj}^n + \frac{r}{6}u_{N,j + 1}^n + \frac{r}{6}u_{N,j - 1}^n} \right]^T},j = 1,2, \ldots ,N - 1.
\end{array}
\]
\[
\begin{array}{l}
u_{h0}^{n + 1} = {\left[ {u_{10}^{n + 1},u_{20}^{n + 1}, \ldots ,u_{N - 1,0}^{n + 1}} \right]^T},\;u_{hN}^{n + 1} = {\left[ {u_{1N}^{n + 1},u_{2N}^{n + 1}, \ldots ,u_{N - 1,N}^{n + 1}} \right]^T},\\
u_{h0}^n = {\left[ {u_{10}^n,u_{20}^n, \ldots ,u_{N - 1,0}^n} \right]^T},\;u_{hN}^n = {\left[ {u_{1N}^n,u_{2N}^n, \ldots ,u_{N - 1,N}^n} \right]^T}.
\end{array}\]
and
$${A_1} = \left( {\begin{array}{*{20}{c}}
{\frac{{2 + 10r}}{3}}&{\frac{{1 - 8r}}{{12}}}&{}&{}&{}\\
{\frac{{1 - 8r}}{{12}}}&{\frac{{2 + 10r}}{3}}&{\frac{{1 - 8r}}{{12}}}&{}&{}\\
{}&{\frac{{1 - 8r}}{{12}}}& \ddots & \ddots &{}\\
{}&{}& \ddots & \ddots &{\frac{{1 - 8r}}{{12}}}\\
{}&{}&{}&{\frac{{1 - 8r}}{{12}}}&{\frac{{2 + 10r}}{3}}
\end{array}} \right),{A_2} = \left( {\begin{array}{*{20}{c}}
{\frac{{2 - 10r}}{3}}&{\frac{{1 + 8r}}{{12}}}&{}&{}&{}\\
{\frac{{1 + 8r}}{{12}}}&{\frac{{2 - 10r}}{3}}&{\frac{{1 + 8r}}{{12}}}&{}&{}\\
{}&{\frac{{1 + 8r}}{{12}}}& \ddots & \ddots &{}\\
{}&{}& \ddots & \ddots &{\frac{{1 + 8r}}{{12}}}\\
{}&{}&{}&{\frac{{1 + 8r}}{{12}}}&{\frac{{2 - 10r}}{3}}
\end{array}} \right);$$
$${B_1} = \left( {\begin{array}{*{20}{c}}
{\frac{{8r - 1}}{{12}}}&{\frac{r}{6}}&{}&{}&{}\\
{\frac{r}{6}}&{\frac{{8r - 1}}{{12}}}&{\frac{r}{6}}&{}&{}\\
{}&{\frac{r}{6}}& \ddots & \ddots &{}\\
{}&{}& \ddots & \ddots &{\frac{r}{6}}\\
{}&{}&{}&{\frac{r}{6}}&{\frac{{8r - 1}}{{12}}}
\end{array}} \right),\;\;\;{B_2} = \left( {\begin{array}{*{20}{c}}
{\frac{{8r + 1}}{{12}}}&{\frac{r}{6}}&{}&{}&{}\\
{\frac{r}{6}}&{\frac{{8r + 1}}{{12}}}&{\frac{r}{6}}&{}&{}\\
{}&{\frac{r}{6}}& \ddots & \ddots &{}\\
{}&{}& \ddots & \ddots &{\frac{r}{6}}\\
{}&{}&{}&{\frac{r}{6}}&{\frac{{8r + 1}}{{12}}}
\end{array}} \right).$$
This is the compact difference scheme (H-OCD) on the equation (\ref{eq3.1}). Similar to Ref. \cite{Tian},
one can prove that the truncation errors of this compact difference method are $O(\tau^2+h^4)$.

\subsection{Two-dimensional Numerical Gradient Scheme}

\quad Next, analogous to the Section 2, let us consider the two-dimensional numerical gradient scheme on the above discrete form (\ref{eq3.15})
by the local Hermite interpolation and collocation polynomial. The intermediate points
$u(x_{i+\frac{1}{2}},y_{j+\frac{1}{2}})$ can be expressed (see (\ref{eq3.20})) by the
values of the mesh points and their partial derivatives $P(x_i,y_j)$ (i.e., $K_{ij}$ and $L_{ij}$, see (\ref{eq3.26}) and (\ref{eq3.27})) around it (see Fig. 1), where $P(x_i,y_j)$ is computed by the difference points around
$u(x_i,y_j)$ (see Fig. 2).

\begin{figure}[htbp]\label{fig6}
\centerline{\includegraphics[width=10.00in,height=3.00in]{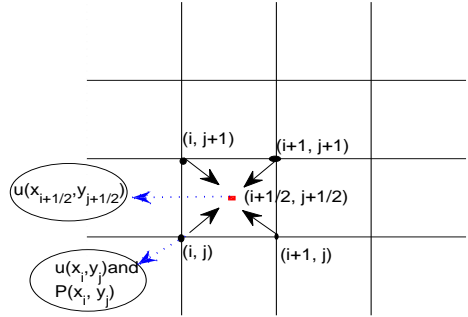}}
\caption{ \small The relationships between
$u(x_{i+\frac{1}{2}},y_{j+\frac{1}{2}})$ and the around points $u(x_i,y_j)$, the
partial derivatives $P(x_i,y_j)$. }
\end{figure}

\begin{figure}[htbp]\label{fig7}
\centerline{\includegraphics[width=10.00in,height=3.00in]{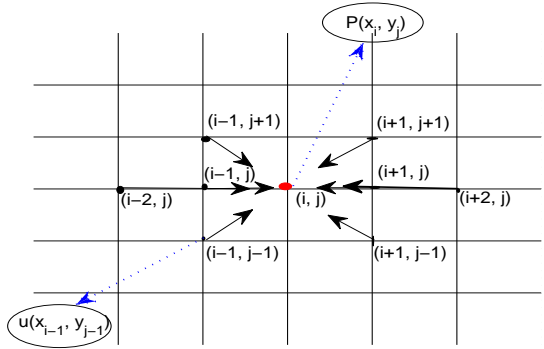}}
\caption{\small The relationships between $P(x_i,y_j)$ and the around points $u(x_i,y_j)$.}
\end{figure}

\subsubsection{The Local Hermite Interpolation and Refinement in two-dimensional case}

\quad For convenience, we denote
$$K_{i,j} = \frac{{\partial u\left( {x_i,y_j,t} \right)}}{{\partial x}}{{{\left. { \buildrel \Delta \over = \frac{{\partial u}}{{\partial x}}} \right|}_{\left( {{x_i},{y_j},t} \right)}}},\;\;L_{i,j} = \frac{{\partial u\left( {x_i,y_j,t} \right)}}{{\partial y}}.$$

Let us consider Hermite bilinear interpolation functions $\Psi _H(x,y,t)$ on the rectangular mesh
${[{x_i},{x_{i + 1}]}}\times {[{y_j},{y_{j + 1}]}}\subset\Omega_h$, and its four vertexes are as follows.
\[
{z_1}: \left( {{x_i},{y_j},t} \right),{z_2}: \left( {{x_{i + 1}},{y_j},t} \right),{z_3}: \left( {{x_i},{y_{j + 1}},t} \right),{z_4}: \left( {{x_{i + 1}},{y_{j + 1}},t} \right) \in {\Omega _h}.
\]
On the segment $z_1-z_2$, let the bilinear interpolation function satisfy the following conditions
\[
    \Psi _H(z_1)=u(z_1), \; \Psi _H(z_2)=u(z_2),
\]
\[
    (\Psi _H)_x(z_1)=u_x(z_1),\; (\Psi _H)_x(z_2)=u(z_2).
\]
Based on Ref.\cite{T.S}, we can get the following Hermite interpolation polynomial
\begin{equation}\label{eq3.19}
      \Psi_H(x_{i+\frac{1}{2}},y_j,t)=\frac{{1}}{{2}}[u(z_1)+u(z_2)]+\frac{{h}}{{8}}[u_x(z_1)-u_x(z_2)],
\end{equation}
where $i=1,2,...,N-1$. The interpolation errors are
\begin{equation}\label{eq3.20}
   \begin{array}{lll}
{\Psi _H}\left( {{x_{i + \frac{1}{2}}},{y_j},t} \right) - u\left( {{x_{i + \frac{1}{2}}},{y_j},t} \right) &= &\frac{{{u_{xxxx}}\left( {{\xi _i},{y_j},t} \right)}}{{4!}}{\left( {{x_{i + \frac{1}{2}}} - {x_i}} \right)^2}{\left( {{x_{i + \frac{1}{2}}} - {x_{i + 1}}} \right)^2}\\
& =& \frac{{h_x^4}}{{384}}{u_{xxxx}}\left( {{\xi _i},{y_j},t} \right), i,j = 1,2, \ldots ,N - 1.
   \end{array}
\end{equation}
where $\xi_j$ lies between $z_1$ and $z_2$ (see, \cite{T.S}).
Thus, we obtain the following approximate computation formula for any $i = 2,3, \ldots ,,N - 1,j = 1,2, \ldots ,N - 1$.
\begin{equation}\label{eq3.21}
{u_{i + \frac{1}{2},j}^k} = \frac{1}{2}\left( {u_{ij}^k + u_{i + 1,j}^k} \right) + \frac{{{h_x}}}{8}\left( {{K_{ij}} - {K_{i + 1,j}}} \right).
\end{equation}
Similarly, we have also that
\begin{equation}\label{eq3.22}
{u_{i,j+\frac{1}{2}}^k} = \frac{1}{2}\left( {u_{ij}^k + u_{i,j + 1}^k} \right) + \frac{{{h_y}}}{8}\left( {{L_{ij}} - {L_{i,j + 1}}} \right).
\end{equation}
Therefore, for $i,j = 2,3, \ldots ,N - 3$, $u_{i + \frac{1}{2},j + \frac{1}{2}}^k $ can be approximated as follows
\begin{equation}\label{eq3.23}
\begin{array}{l}
u_{i + \frac{1}{2},j + \frac{1}{2}}^k = \frac{1}{2}\left( {u_{i,j + \frac{1}{2}}^k + u_{i + 1,j + \frac{1}{2}}^k} \right) + \frac{h}{8}\left( {{K_{i,j + \frac{1}{2}}} - {K_{i + 1,j + \frac{1}{2}}}} \right)\\
 = \frac{1}{4}\left( {u_{ij}^k + u_{i,j + 1}^k} \right) + \frac{h}{{16}}\left( {{L_{ij}} - {L_{i,j + 1}} + {L_{i + 1,j}} - {L_{i + 1,j + 1}}} \right) + \\
\frac{1}{4}\left( {u_{i + 1,j}^k + u_{i + 1,j + 1}^k} \right) + \frac{h}{{16}}\left( {{K_{ij}} - {K_{i + 1j}} + {K_{ij + 1}} - {K_{i + 1j + 1}}} \right).
\end{array}
\end{equation}

In Section 3.2, we will prove that the above refinement scheme has the four-order accuracy in space direction, see Theorem 3.2.

\subsubsection{The Collocation Polynomial in two-dimensional case}

\quad Next, we use the collocation polynomial
method to obtain the approximate values of $K_{ij}$ and $L_{ij}$.
For convenience, we consider the sub-domain
\[
   [x_{i-1},x_{i+1}]\times [y_{j-1},y_{j+1}]\subset\Omega, \; i,j=1,2,...,N-1.
\]
and denote
\[
\xi  = x - {x_i},\eta  = y - {y_j},\left( {x,y} \right) \in {\Omega _h},\left( {i,j = 1,2, \ldots ,N - 1} \right).
\]
In order to get the approximation polynomial of $u$, we consider the polynomial space
\begin{equation}
   {H_4} = span\left\{ {1,\xi ,\eta ,{\xi ^2},\xi \eta ,{\eta ^2},{\xi ^3},{\xi ^2}\eta ,\xi {\eta ^2},{\eta ^3},{\xi ^4},{\xi ^2}{\eta ^2},{\eta ^4}} \right\}
\end{equation}
and define the approximation polynomial as follows
\begin{equation}
   \begin{array}{l}
H(\xi ,\eta ) = {a_0} + {a_1}\xi  + {a_2}\eta  + {a_3}{\xi ^2} + {a_4}\xi \eta  + {a_5}{\eta ^2} + {a_6}{\xi ^3} + {a_7}{\xi ^2}\eta  + \\
\quad\quad\quad \quad \;\;{a_8}\xi {\eta ^2} + {a_9}{\eta ^3} + {a_{10}}{\xi ^4} + {a_{11}}{\xi ^2}{\eta ^2} + {a_{12}}{\eta ^4}.
\end{array}
\end{equation}

\quad Let
\[\begin{array}{l}
H\left( {{x_m},{y_n}} \right) = u\left( {{x_m},{y_n}} \right),m = i - 1,i,i + 1,n = j - 1,j,j + 1,{\kern 1pt} \left( {m - i} \right)\left( {n - j} \right) \ne 0;\\
\frac{{{\partial ^2}H\left( {{x_m},{y_n}} \right)}}{{\partial {x^2}}} + \frac{{{\partial ^2}H\left( {{x_m},{y_n}} \right)}}{{\partial {y^2}}} = \frac{{\partial u}}{{\partial t}}\left( {{x_m},{y_n},t} \right),
\end{array}\]
where $n = j,m = i - 1,i,i + 1$ and $n = j \pm 1,m = i$. Then, we can obtain the following numerical gradient approximate scheme.
\begin{equation}\label{eq3.26}
\begin{array}{l}
{K_{ij}} = \frac{1}{{12h}}\left[ {8u\left( {{x_{i + 1}},{y_j},t} \right) - 8u\left( {{x_{i - 1}},{y_j},t} \right) + u\left( {{x_{i - 2}},{y_j},t} \right) - u\left( {{x_{i + 2}},{y_j},t} \right)} \right],\\
\quad\quad \quad i = 2,3, \ldots ,N - 2,j = 1,2, \ldots ,N - 1;
\end{array}
\end{equation}

\begin{equation}\label{eq3.27}
\begin{array}{l}
{L_{ij}} = \frac{1}{{12h}}\left[ {8u\left( {{x_i},{y_{j + 1}},t} \right) - 8u\left( {{x_i},{y_{j - 1}},t} \right) + u\left( {{x_i},{y_{j - 2}},t} \right) - u\left( {{x_i},{y_{j + 2}},t} \right)} \right],\\
\quad\quad \quad i = 1,2, \ldots ,N - 1,j = 2,3, \ldots ,N - 2.
\end{array}
\end{equation}

\subsection{The truncation errors of numerical gradient scheme}

\quad As stated in the previous Section 2, the truncation errors of the compact difference method in \cite{Sun.Zhang} are $O(\tau^2+h^4)$.
In fact, the above numerical gradient schemes (\ref{eq3.26}) and (\ref{eq3.27}) have also the same convergence order.

\begin{theorem}
If $ u(x,y,t)\in \mathbb{C}^6_x(\Omega) $ and $ u(x,y,t)\in \mathbb{C}^6_y(\Omega) $, then
we have
\begin{equation}
\left| {{K_{ij}} - \frac{{\partial u\left( {{x_i},{y_j},t} \right)}}{{\partial x}}} \right| < O\left( {{h^4}} \right),(i = 2,3, \ldots ,N - 2,j = 1,2, \ldots ,N - 1),
\end{equation}

\begin{equation}
\left| {{L_{ij}} - \frac{{\partial u\left( {{x_i},{y_j},t} \right)}}{{\partial y}}} \right| < O\left( {{h^4}} \right),\left( {j = 2,3, \ldots ,N - 2,i = 1,2, \ldots ,N - 1} \right).
\end{equation}
\end{theorem}

\textbf{Proof.} According to (\ref{eq3.26}), we know
\[
   \begin{array}{lll}
     {K_{ij}} &=& \frac{1}{{12h}}\left[ {8u\left( {{x_{i + 1}},{y_j},t} \right) - 8u\left( {{x_{i - 1}},{y_j},t} \right) + u\left( {{x_{i - 2}},{y_j},t} \right) - u\left( {{x_{i + 2}},{y_j},t} \right)} \right]\\
     &=&\frac{4}{3}\frac{{u\left( {{x_{i + 1}},{y_j},t} \right) - u\left( {{x_{i - 1}},{y_j},t} \right)}}{{2h}} - \frac{1}{3}\frac{{u\left( {{x_{i - 2}},{y_j},t} \right) - u\left( {{x_{i + 2}},{y_j},t} \right)}}{{4h}}.
   \end{array}
\]
In addition, according to Taylor series expansion theorem, we have
\[{\frac{{\partial u\left( {{x_i},{y_j},t} \right)}}{{\partial x}} = \frac{4}{3}\frac{{u\left( {{x_{i + 1}},{y_j},t} \right) - u\left( {{x_{i - 1}},{y_j},t} \right)}}{{2h}} - \frac{1}{3}\frac{{u\left( {{x_{i - 2}},{y_j},t} \right) - u\left( {{x_{i + 2}},{y_j},t} \right)}}{{4h}} + O\left( {{h^4}} \right)}\]
So
\[\left| {\frac{{\partial u\left( {{x_i},{y_j},t} \right)}}{{\partial x}} - K_{ij}} \right| = O\left( {{h^4}} \right).
\]
Similarly, we may also prove that
\[\left| {\frac{{\partial u\left( {{x_i},{y_j},t} \right)}}{{\partial y}} - L_{ij}} \right| = O\left( {{h^4}} \right).
\]
Thus the proof is completed. $\Box$

\quad By the above theorem, we know that the accuracy of numerical gradient schemes (\ref{eq3.26}) and (\ref{eq3.27}) is $O(h^4)$
in the space direction. In fact, in the intermediate points $({x_{i + \frac{1}{2}}},{y_{j + \frac{1}{2}}},t)$, the above refinement scheme (\ref{eq3.23}) has also the four-order accuracy in space direction.

\begin{theorem}
   If $ u(x,y,t)\in \mathbb{C}^6_x(\Omega) $ and $ u(x,y,t)\in \mathbb{C}^6_y(\Omega) $,  then
 \begin{equation}
 \left| {u\left( {{x_{i + \frac{1}{2}}},{y_{j + \frac{1}{2}}},{t_k}} \right) - u_{i + \frac{1}{2},j + \frac{1}{2}}^k} \right| \le O \left( {{h^4}} \right),i,j = 2,3, \ldots ,N - 3.
 \end{equation}
\end{theorem}

\textbf{Proof.} First, by the Taylor expansion of $u\left( {{x_i},{y_j},{t_k}} \right)$ at $(x_{j+\frac{1}{2}},T)$, we have
\[
\begin{array}{l}
u\left( {{x_i},{y_j},{t_k}} \right) = u\left( {{x_{i + \frac{1}{2}}},{y_j},{t_k}} \right) - \frac{h}{2}{u_x}\left( {{x_{i + \frac{1}{2}}},{y_j},{t_k}} \right) + \frac{{{h^2}}}{8}{u_{xx}}\left( {{x_{i + \frac{1}{2}}},{y_j},{t_k}} \right) - \\
{\kern 1pt} {\kern 1pt} {\kern 1pt} {\kern 1pt} {\kern 1pt} {\kern 1pt} {\kern 1pt} {\kern 1pt} {\kern 1pt} {\kern 1pt} {\kern 1pt} {\kern 1pt} {\kern 1pt} {\kern 1pt} {\kern 1pt} {\kern 1pt} {\kern 1pt} {\kern 1pt} {\kern 1pt} {\kern 1pt} {\kern 1pt} {\kern 1pt} {\kern 1pt} {\kern 1pt} {\kern 1pt} {\kern 1pt} {\kern 1pt} {\kern 1pt} {\kern 1pt} {\kern 1pt} {\kern 1pt} {\kern 1pt} {\kern 1pt} {\kern 1pt} {\kern 1pt} {\kern 1pt} {\kern 1pt} {\kern 1pt} {\kern 1pt} {\kern 1pt} {\kern 1pt} {\kern 1pt} {\kern 1pt} {\kern 1pt} {\kern 1pt} {\kern 1pt} {\kern 1pt} {\kern 1pt} {\kern 1pt} {\kern 1pt} {\kern 1pt} {\kern 1pt} {\kern 1pt} {\kern 1pt} {\kern 1pt} {\kern 1pt} {\kern 1pt} {\kern 1pt} {\kern 1pt} {\kern 1pt} {\kern 1pt} {\kern 1pt} {\kern 1pt} {\kern 1pt} {\kern 1pt} {\kern 1pt} {\kern 1pt} {\kern 1pt} {\kern 1pt} {\kern 1pt} {\kern 1pt} \frac{{{h^3}}}{{48}}{u_{xxx}}\left( {{x_{i + \frac{1}{2}}},{y_j},{t_k}} \right) + O\left( {{h^4}} \right),
\end{array}
\]

\[\begin{array}{l}
u\left( {{x_{i + 1}},{y_j},{t_k}} \right) = u\left( {{x_{i + \frac{1}{2}}},{y_j},{t_k}} \right) + \frac{h}{2}{u_x}\left( {{x_{i + \frac{1}{2}}},{y_j},{t_k}} \right) + \frac{{{h^2}}}{8}{u_{xx}}\left( {{x_{i + \frac{1}{2}}},{y_j},{t_k}} \right) + \\
{\kern 1pt} {\kern 1pt} {\kern 1pt} {\kern 1pt} {\kern 1pt} {\kern 1pt} {\kern 1pt} {\kern 1pt} {\kern 1pt} {\kern 1pt} {\kern 1pt} {\kern 1pt} {\kern 1pt} {\kern 1pt} {\kern 1pt} {\kern 1pt} {\kern 1pt} {\kern 1pt} {\kern 1pt} {\kern 1pt} {\kern 1pt} {\kern 1pt} {\kern 1pt} {\kern 1pt} {\kern 1pt} {\kern 1pt} {\kern 1pt} {\kern 1pt} {\kern 1pt} {\kern 1pt} {\kern 1pt} {\kern 1pt} {\kern 1pt} {\kern 1pt} {\kern 1pt} {\kern 1pt} {\kern 1pt} {\kern 1pt} {\kern 1pt} {\kern 1pt} {\kern 1pt} {\kern 1pt} {\kern 1pt} {\kern 1pt} {\kern 1pt} {\kern 1pt} {\kern 1pt} {\kern 1pt} {\kern 1pt} {\kern 1pt} {\kern 1pt} {\kern 1pt} {\kern 1pt} {\kern 1pt} {\kern 1pt} {\kern 1pt} {\kern 1pt} {\kern 1pt} {\kern 1pt} {\kern 1pt} {\kern 1pt} {\kern 1pt} {\kern 1pt} {\kern 1pt} {\kern 1pt} {\kern 1pt} {\kern 1pt} {\kern 1pt} {\kern 1pt} {\kern 1pt} {\kern 1pt} \frac{{{h^3}}}{{48}}{u_{xxx}}\left( {{x_{i + \frac{1}{2}}},{y_j},{t_k}} \right) + O\left( {{h^4}} \right).
\end{array}
\]
Therefore,
\[
\begin{array}{l}
u\left( {{x_{i + \frac{1}{2}}},{y_{j + \frac{1}{2}}},{t_k}} \right) = \frac{1}{4}\left[ {u\left( {{x_i},{y_j},{t_k}} \right) + u\left( {{x_i},{y_{j + 1}},{t_k}} \right) + u\left( {{x_{i + 1}},{y_j},{t_k}} \right) + u\left( {{x_{i + 1}},{y_{j + 1}},{t_k}} \right)} \right] - \\
\;\; \;\;\frac{{{h^2}}}{8}{u_{xx}}\left( {{x_{i + \frac{1}{2}}},{y_{j + \frac{1}{2}}},{t_k}} \right) - \frac{{{h^2}}}{{16}}[{u_{yy}}\left( {{x_i},{y_{j + \frac{1}{2}}},{t_k}} \right) - {u_{yy}}\left( {{x_{i + 1}},{y_{j + \frac{1}{2}}},{t_k}} \right)] + O\left( {{h^4}} \right)\\
 = \frac{1}{4}\left[ {u\left( {{x_i},{y_j},{t_k}} \right) + u\left( {{x_i},{y_{j + 1}},{t_k}} \right) + u\left( {{x_{i + 1}},{y_j},{t_k}} \right) + u\left( {{x_{i + 1}},{y_{j + 1}},{t_k}} \right)} \right] - \\
\;\; \;\;\frac{h}{{16}}\frac{\partial }{{\partial y}}{\delta _y}\left[ {u\left( {{x_i},{y_{j + \frac{1}{2}}},{t_k}} \right) - u\left( {{x_{i + 1}},{y_{j + \frac{1}{2}}},{t_k}} \right)} \right] - \frac{h}{8}\frac{\partial }{{\partial x}}\left[ {{\delta _x}u\left( {{x_{i + \frac{1}{2}}},{y_{j + \frac{1}{2}}},{t_k}} \right)} \right] + O\left( {{h^4}} \right)\\
 = \frac{1}{4}\left( {u_{ij}^k + u_{i,j + 1}^k + u_{i + 1,j}^k + u_{i + 1,j + 1}^k} \right) - \frac{h}{{16}}\left( {{L_{ij}} - {L_{i,j + 1}} + {L_{i + 1j}} - {L_{i + 1j + 1}}} \right) - \\
\;\;\;\;\frac{h}{{16}}\left( {{K_{ij}} - {K_{i + 1,j}} + {K_{i,j + 1}} - {K_{i + 1,j + 1}}} \right) + O\left( {{h^4}} \right)\\
 = u_{i + \frac{1}{2},j + \frac{1}{2}}^k + O\left( {{h^4}} \right).
\end{array}
\]
That is, the conclusion holds. $\Box$

\quad  In addition, to reduce the total computing time, we also consider the Richardson extrapolation
on the H-OCD scheme (\ref{eq3.15}) in two-dimensional case. For convenience,
we take the following initial-boundary problem as a simple example.
\begin{equation}\label{eq3.31}
     \begin{array}{lll}
  {u_t} - \Delta u = {F_u}\left( {x,y,t} \right),\left( {x,y,t} \right) \in \left( {a,b} \right) \times \left( {c,d} \right) \times \left( {0,T} \right],\\
u\left( {x,c,t} \right) = u\left( {x,d,t} \right) = u\left( {a,y,t} \right) = u\left( {b,y,t} \right) = 0,\left( {x,y} \right) \in \left[ {a,b} \right] \times \left[ {c,d} \right],0 \le t \le T,\\
u\left( {x,y,0} \right) = 0,\left[ {x,y} \right] \in \left[ {a,b} \right] \times \left[ {c,d} \right]
      \end{array}
\end{equation}
with the smooth solution $u(x,y,t)$, where
\[{F_u}\left( {x,y,t} \right) = \frac{1}{{24}}\frac{{{\partial ^3}u\left( {x,y,t} \right)}}{{\partial {t^3}}} - \frac{1}{8}\frac{{{\partial ^4}u\left( {x,y,t} \right)}}{{\partial {x^2}\partial {t^2}}}.
\]

\begin{theorem}
Let $u\left( {x,y,t} \right) \in {\mathbb{C}^{8,6}}\left( {{\Omega _h} \times \left[ {0,T} \right]} \right)$ be the solution of Eq. (\ref{eq3.1}) with initial-boundary (\ref{eq3.31}) and $u_{ij}^k\left( {h,\tau } \right)$ is the numerical solution of H-OCD scheme (\ref{eq3.15}) with the time step $\tau$ and the space step $h$. Then
\[\mathop {\max }\limits_{1 \le i,j \le N - 1,1 \le k \le M} \left| {u\left( {{x_i},{y_j},{t_k}} \right) - \left[ {\frac{4}{3}u_{ij}^{2k}\left( {h,\frac{\tau }{2}} \right) - \frac{1}{3}u_j^k\left( {h,\tau } \right)} \right]} \right| = O\left( {{\tau ^4} + {h^4}} \right).
\]
\end{theorem}

\textbf{Proof.}  The proof is completely similar to the Theorem 2.3. In addition, corresponding numerical experiments will be shown in Table 10. $\Box$

\section{Numerical Experiments}

\subsection{Numerical Experiments for the One-dimensional Case}
\quad \textbf{Example 4.1} When $u(x,0)=\sin(\pi x),
u(0,t)=u(1,t)=0$ for the equation (\ref{eq1}) with $(x,t)\in
(0,1)\times(0,T]$, the exact solution of the problem (\ref{eq1}) is
\[
   u(x,t)=\exp(-\pi^2t)\sin(\pi x).
\]

\quad Next, let us observe and compare the numerical solutions in the
same number of points and the time for the above two schemes.

\quad First, we note that matrix computations are based on LAPACK, and optimized
basic linear algebra subroutines (BLAS) on all Matlab platforms, which speeds
up matrix multiplications and the LAPACK routines themselves, according to Matlab user manual.
Therefore, all the numerical experiments were performed in MATLAB
2011b. In addition, for convenience, we denote $Rate(h)=\log_{2}(\frac{{Error(h)}}{{Error(\frac{h}{2})}})$ and
$Error(h)=\max_{x_k=x_0+kh,k=0,1,...N}\{\mid(u(x_k,T)-u^T_k)\mid\}$,
where $u(x_k,T)$ represents the exact solution and $u^T_k$ is the
numerical solution. Let $$Error(P)=\max_{x_k=x_0+kh,k=0,1,...N}\{\mid \frac{{\partial u}}{{\partial x}}(x_j,T)-P^T_k)\mid\}.$$

\quad Table 1 lists the computational results of the mesh-grid points,
intermediate points and $u_x$ with different spacial step sizes
when time step size is fixed as $\tau=1/100000.$ We can see that the
convergence orders in space can reach $O(h^4)$ which is consistent with
the theoretical analysis in this article.

\begin{table}[htbp]
\begin{tabular*}{14cm}{p{25pt}p{64pt}p{25pt}p{64pt}p{25pt}p{64pt}p{25pt}}
\multicolumn{7}{p{410pt}}{\small Table 1}\\
\multicolumn{7}{p{410pt}}{\small Errors and rate of intermediate points and numerical gradient $P$ (\ref{eq13}) in space direction with $\tau=1/100000$.}\\
\hline $h$ &\multicolumn{2}{p{50pt}}{mesh-grid points}&\multicolumn{2}{p{50pt}}{intermediate points} &\multicolumn{2}{p{50pt}}{P(i.e., $\frac{{\partial u}}{{\partial x}}$)}\\
\cline{2-2}\cline{4-4}\cline{6-6}&Error&Rate&Error&Rate&Error&Rate\\
\hline 1/4&8.3491e-007&4.0355&5.4790e-007&3.5729&3.7721e-006&3.7762\\
1/8&5.0915e-008&4.0073&4.6043e-008&3.9694&2.5732e-007&3.9370\\
1/16&3.1660e-009&4.0045&2.9394e-009&3.9952&1.7976e-008&3.9823\\
1/32&1.9725e-010&4.0466&1.8433e-010&4.0474&1.1374e-009&3.9715\\
1/64&1.1936e-011&*&1.1148e-011&*&7.2504e-011&*\\
\hline
\end{tabular*}
\label{tab1}
\end{table}

\begin{figure}\label{fig1}
  \centering
    \includegraphics[width=0.90\textwidth]{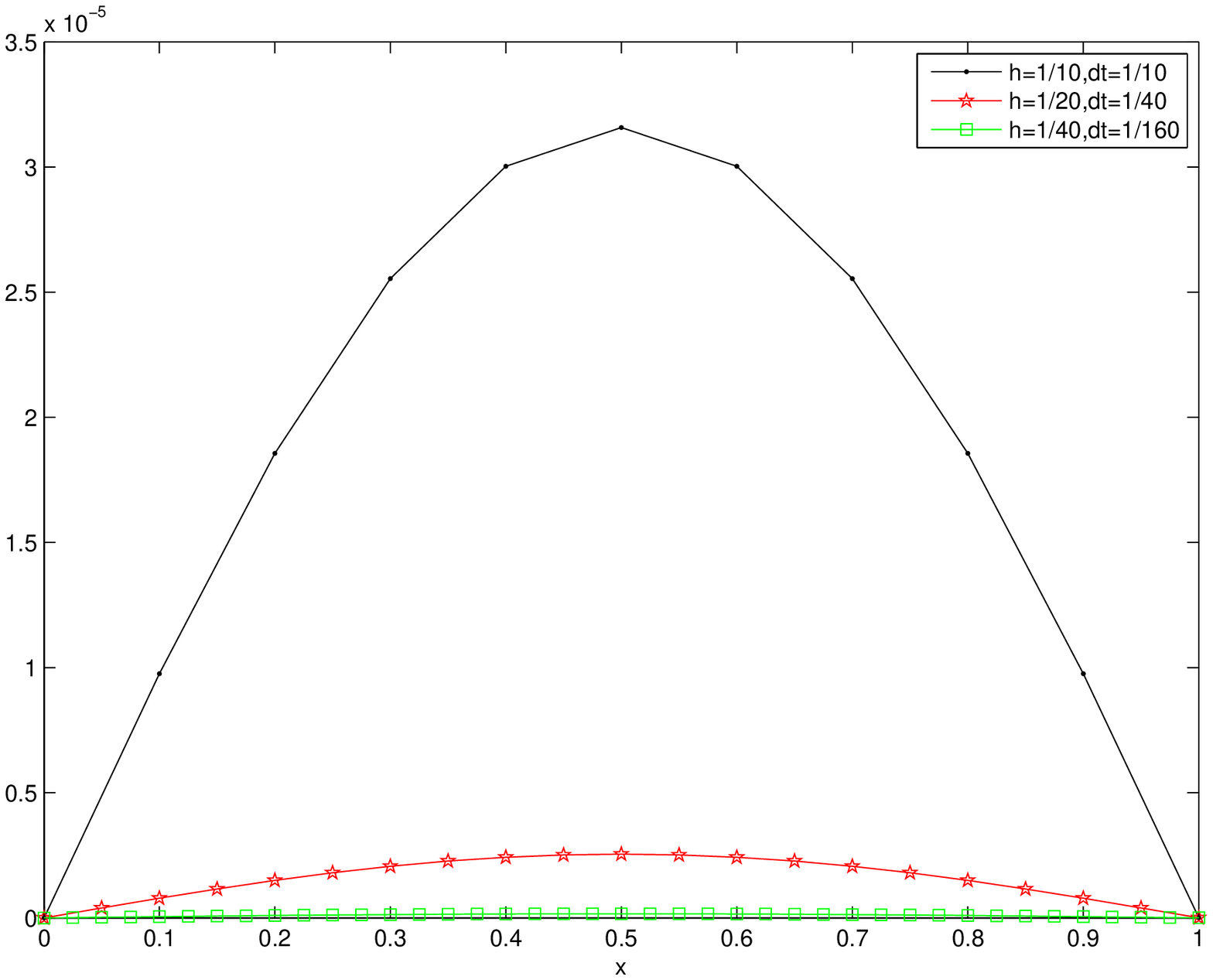}
    \includegraphics[width=0.90\textwidth]{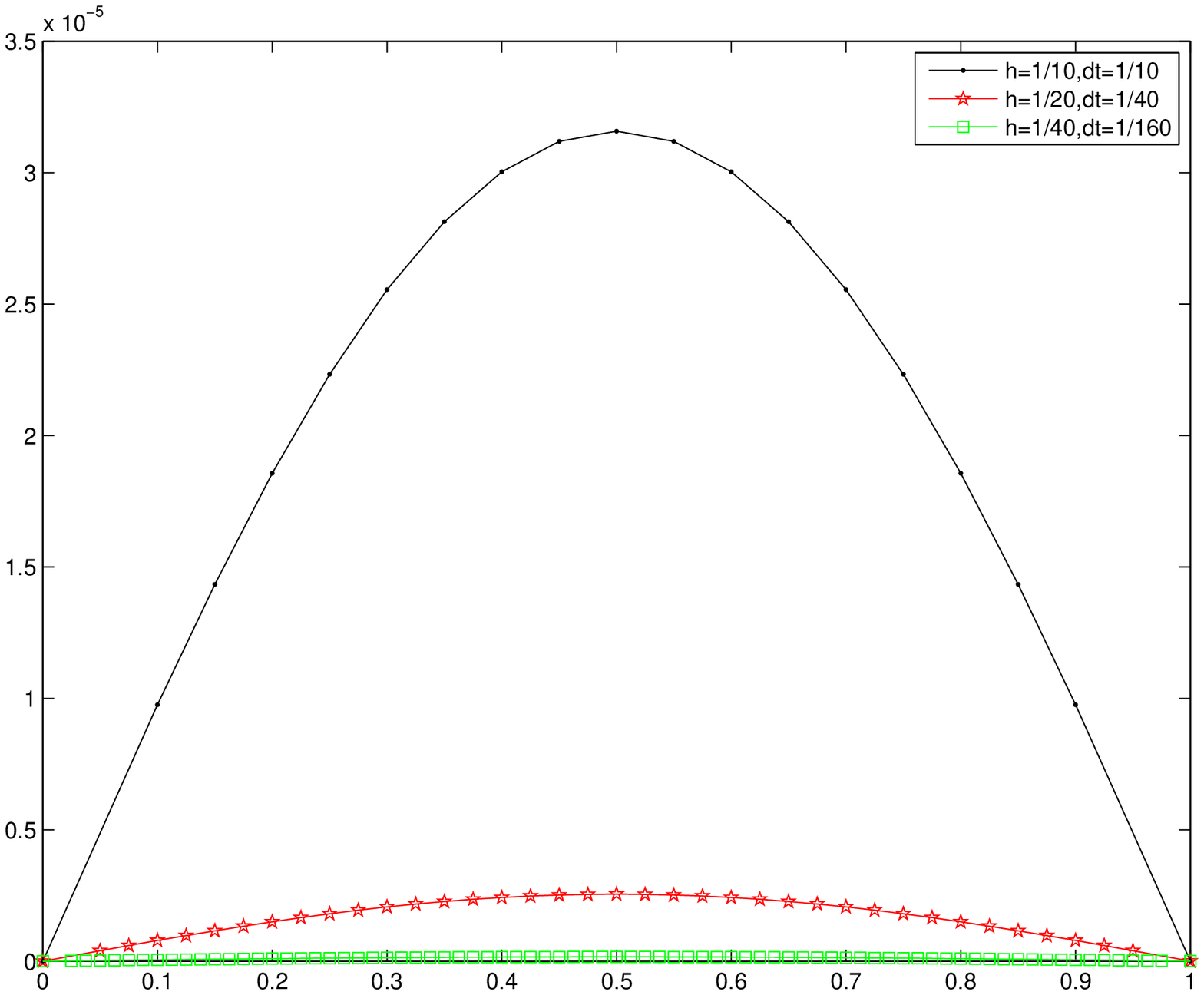}
  \caption{\small Top: The error curves of the mesh-grid points (H-OCD method) when $T=1$;
          Bottom: The error curves of all the points (new method) when $T=1$.}
\end{figure}

\quad Fig.1 displays the errors curves with different step sizes
of the mesh-grid points (by H-OCD method) and all the points (by new method) when $T=1$. They display that,
the changes of the truncation errors in the mesh-grid points and the other
points are all large with the changes of $h$ and $\tau$. At the same time, the
shape of the curves is approximately the same. That means that the points obtained
through the new method are not worse than the H-OCD method.

\begin{figure}[htbp]\label{fig2}
\centerline{\includegraphics[width=8.00in, height=3.90in]{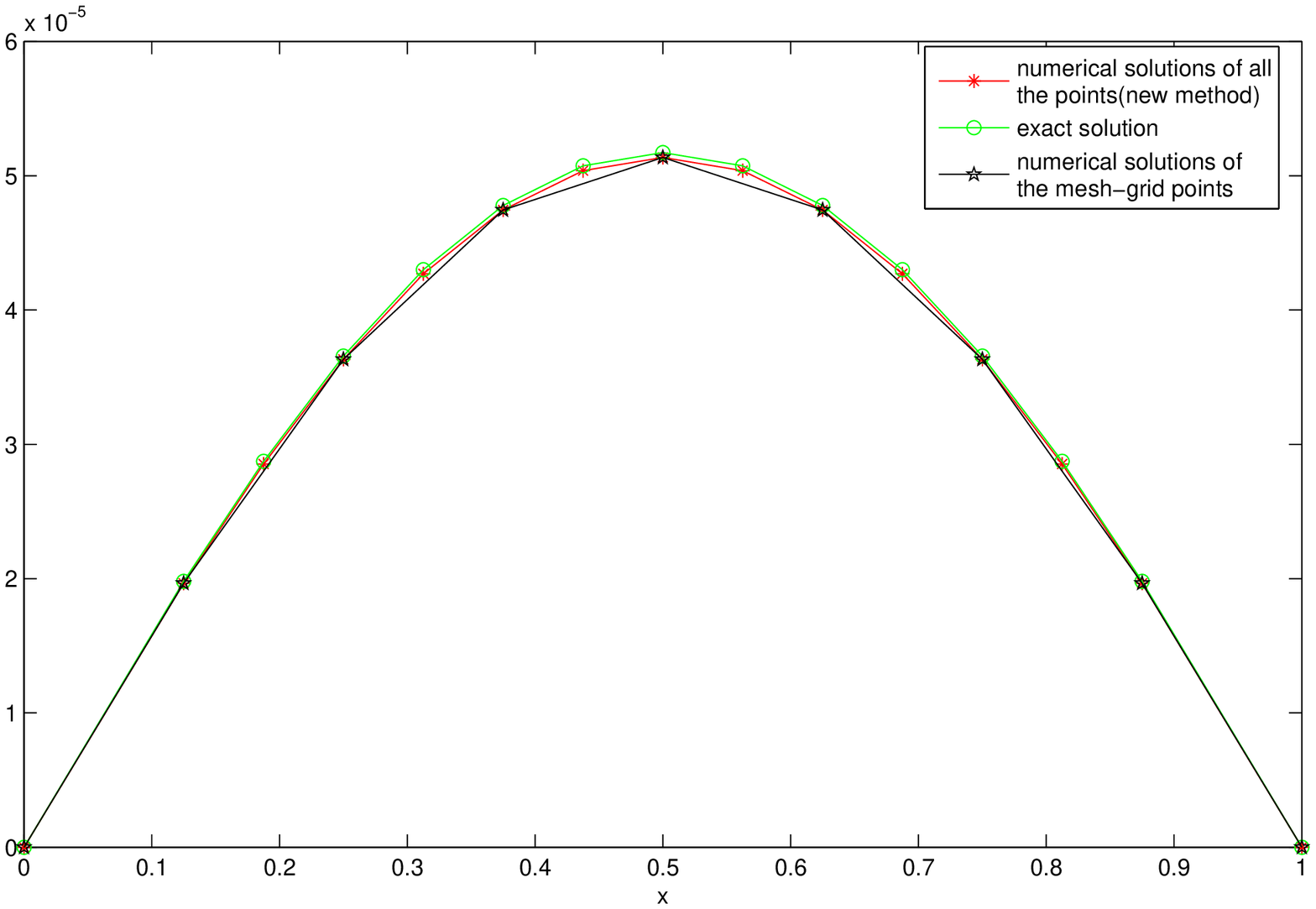}}
\caption{\small The curves of the numerical solutions and exact solutions in all the points (new method)
and the curve of the numerical solutions in the mesh-grid points (H-OCD method) when $h=1/8, \tau=1/100$, $t=1$. }
\end{figure}

\quad Fig.2 shows that the curve of the numerical solutions \textbf{(the red line)} in all the points (by new method)
is \textbf{more close} to the curve of the exact solutions \textbf{(the green line)} when $h=1/8,\ \tau=1/100, T=1$.
That is to say, the simulation result of the red line is better than another. In order to make the
figure 2 more clearly, we choose $h=1/8$.

\begin{figure}\label{fig3}
  \centering
    \includegraphics[width=0.90\textwidth]{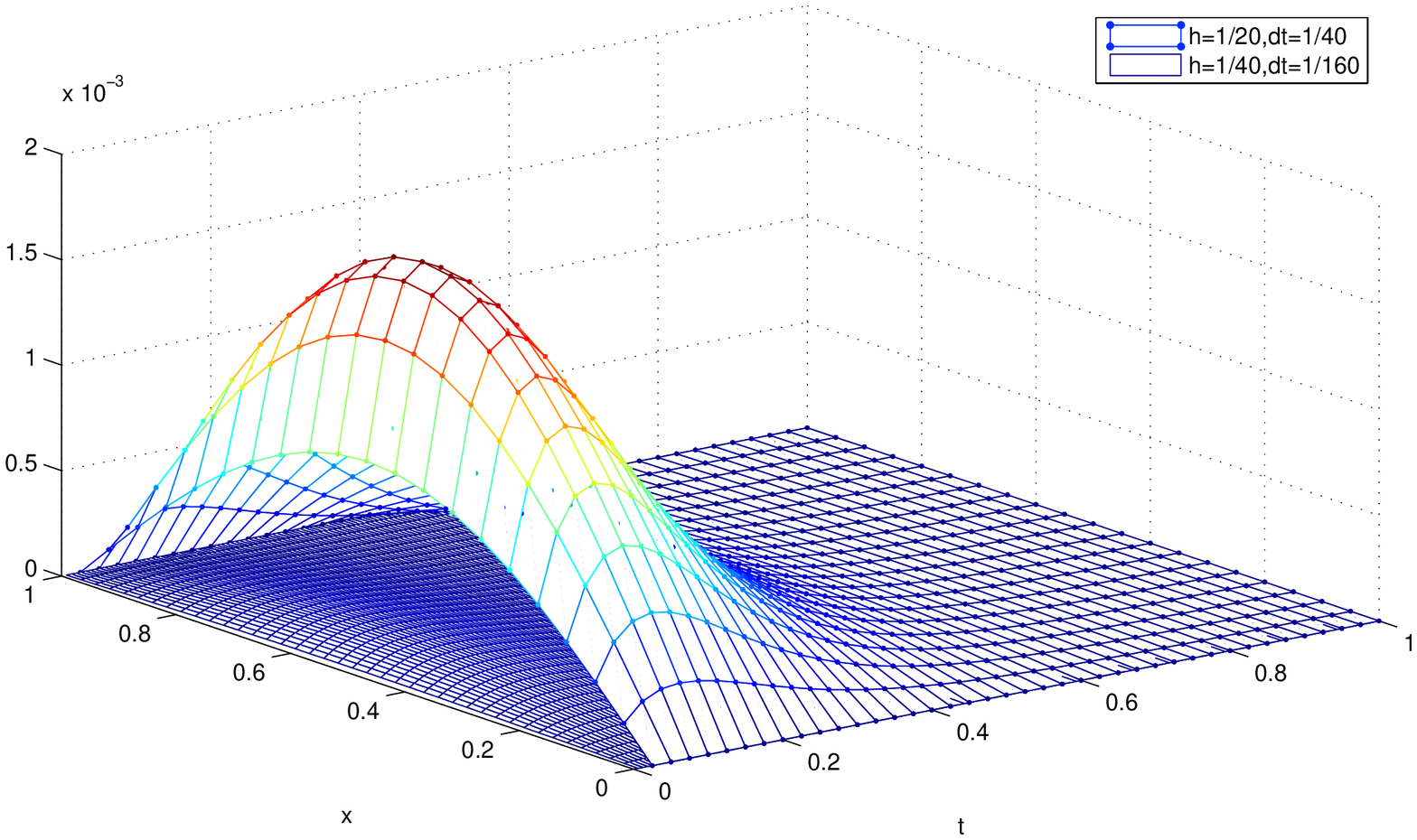}
    \includegraphics[width=0.90\textwidth]{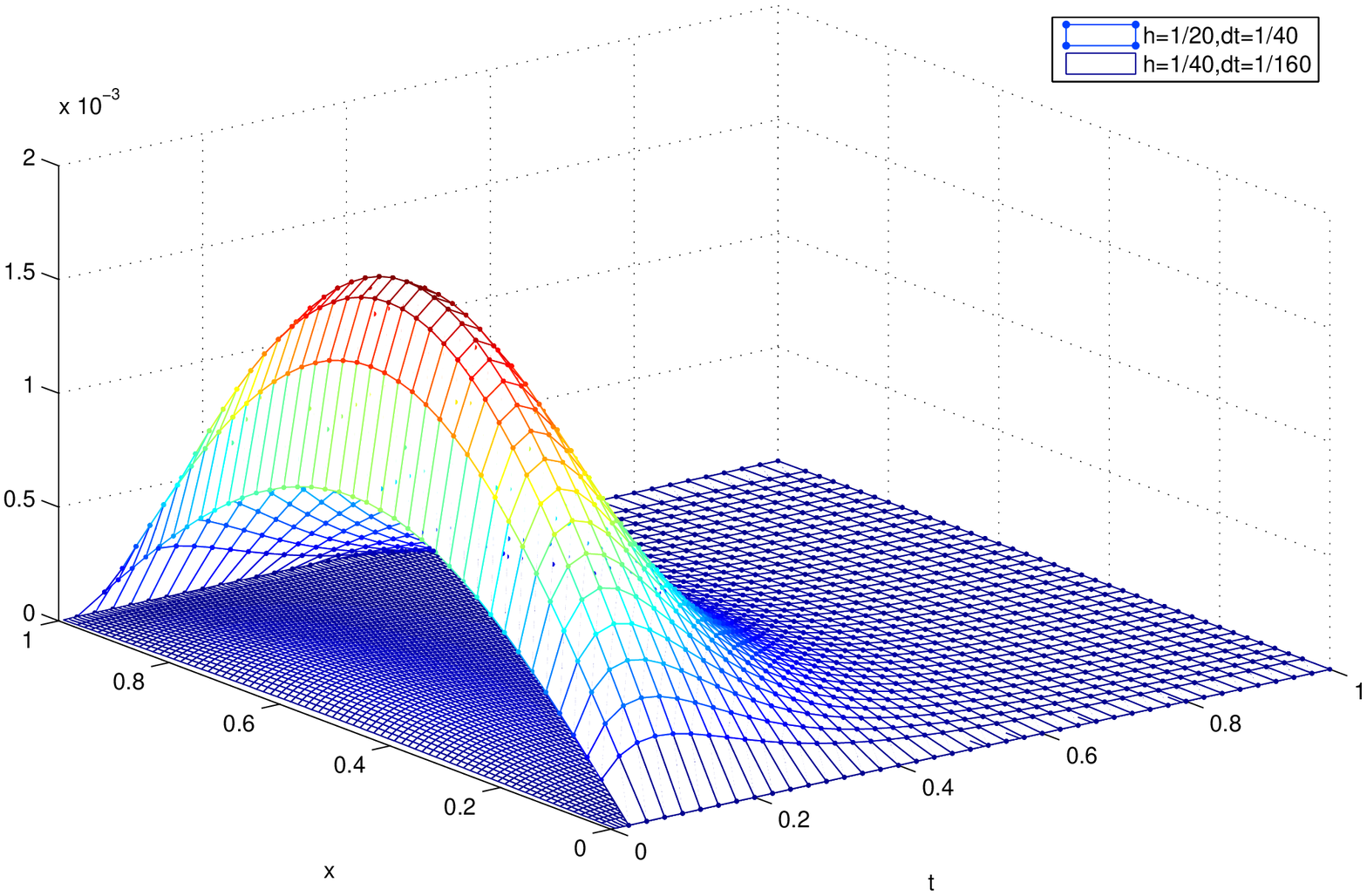}
  \caption{\small Top:  The error surface map of the compact difference scheme in the mesh-grid points;
          Bottom: The error surface map of the numerical gradient scheme in the intermediate points.}
\end{figure}

\quad  Fig.3 displays the error surface maps with different step sizes
in both spatial and time directions of the mesh-grid points (by H-OCD method) and all the
points (by new method) when $t=1$. They display that, the changes of the truncation
errors in the mesh-grid points and the other points are large with the changes
of $h$ and $\tau$. At the same time, the shape of the curves is approximately
the same. That means that the points obtained through the new method are very well too.

\begin{table}[htbp]
\begin{tabular*}{14cm}{p{47pt}p{76pt}p{43pt}p{86pt}p{43pt}}
\multicolumn{5}{p{410pt}}{\small Table 2}\\
\multicolumn{5}{p{410pt}}{\small Errors of the numerical solutions in the mesh-grid points (H-OCD method), all the points (numerical gradient scheme) and the time to get those solutions when $\tau=h^2, T=1, h=1/(n-1)$ .}\\
\hline grid node number&\multicolumn{2}{p{76pt}}{H-OCD method}&\multicolumn{2}{p{86pt}}{numerical gradient scheme}\\
\cline{2-2}\cline{4-4}&Error&Time&Error&Time\\
\hline N=15&6.0041e-008&0.1544&9.5518e-007&0.0312\\
N=31&3.7541e-009&0.5725&6.0041e-008&0.1560\\
N=63&2.3464e-010&2.5389&3.7623e-009& 0.5839\\
N=127&1.4665e-011&11.505&2.3536e-010& 2.7233\\
N=255& 9.1633e-013&\fbox{\small 142.64}&1.4713e-011& \fbox{12.8879}\\
\hline
\end{tabular*}
\label{tab2}
\end{table}

\quad In addition, from Table 2, we can know that, H-OCD method takes more time
to compute the same number of the difference points, compared
to the new method. For example, if we need the numerical solutions
of 255 points to simulate the real figure, we only need $h=1/128$.
Through the method talked in this article, we can get the numerical
solutions of 255 points. However, it just take 12.8879 seconds to do
that, which is much less than 142.64 seconds the H-OCD method needs.

\quad \textbf{Example 4.2} For $u(x,0)=\exp(x), u(0,t)=\exp(t),
u(1,t)=\exp(1+t)$ with $(x,t)\in (0,1)\times(0,T]$, the exact
solution of the problem (\ref{eq1}) is
\[
   u(x,t)=\exp(x+t).
\]

\quad Next, let us compare the numerical solution with the exact
solution as follows (see, Tab. 3-6).

\begin{table}[htbp]
\begin{tabular*}{14cm}{p{30pt}p{60pt}p{30pt}p{60pt}p{30pt}p{60pt}p{30pt}}
\multicolumn{7}{p{410pt}}{\small Table 3}\\
\multicolumn{7}{p{410pt}}{\small Errors and rate of H-OCD scheme (\ref{eq7}-\ref{eq9}), intermediate points (new method) and numerical gradient $P_j$ (\ref{eq13}) in space direction with $\tau=1/100000$.}\\
\hline $h$ &\multicolumn{2}{p{50pt}}{mesh-grid points} &\multicolumn{2}{p{50pt}}{intermediate points}&\multicolumn{2}{p{70pt}}{$P_j$(i.e., $\frac{{\partial u}}{{\partial x}}$)}\\
\cline{2-2}\cline{4-4}\cline{6-6}&Error&Rate&Error&Rate&Error&Rate\\
\hline 1/4&8.4064e-006&3.9974&2.9136e-004&3.8091&6.8324e-003&2.7543\\
1/8&5.2636e-007&3.9895&2.0787e-005&3.9099&1.1027e-003&2.8763\\
1/16&3.3138e-008&3.9993&1.3828e-006&3.9564&1.3792e-004&2.9379\\
1/32&2.0721e-009&4.0472&8.9081e-008&3.9787&1.7998e-005&2.9689\\
1/64&1.2534e-010&*&5.6505e-009&*&2.2987e-006&*\\
\hline
\end{tabular*}
\label{tab3}
\end{table}

\begin{table}[htbp]
\begin{tabular*}{14cm}{p{30pt}p{60pt}p{30pt}p{60pt}p{30pt}p{60pt}p{30pt}}
\multicolumn{7}{p{410pt}}{\small Table 4}\\
\multicolumn{7}{p{410pt}}{\small Errors and rate of H-OCD scheme (\ref{eq7}-\ref{eq9}), intermediate points (new method) and numerical gradient $P_j$ (\ref{eq13}) in time direction with $ h=1/10000$.}\\
\hline $\tau$ &\multicolumn{2}{p{50pt}}{Compact difference} &\multicolumn{2}{p{50pt}}{intermediate points}&\multicolumn{2}{p{70pt}}{$P_j$(i.e., $\frac{{\partial u}}{{\partial x}}$)}\\
\cline{2-2}\cline{4-4}\cline{6-6}&Error&Rate&Error&Rate&Error&Rate \\
\hline 1/10&4.3449e-004&1.9988&4.3449e-004&1.9988&2.0491e-003&1.9769\\
1/20&1.0871e-004&1.9998&1.0871e-004&1.9998&5.2055e-004&1.9887\\
1/40&2.7183e-005&1.9999&2.7183e-005&1.9999&1.3116e-004&1.9945\\
1/80&6.7960e-006&2.0005&6.7960e-006&2.0005&3.2914e-005&1.9987\\
1/160&1.6984e-006&2.0053&1.6984e-006&2.0053&8.2362e-006&2.0053\\
1/320&4.2303e-007&2.0246&4.2303e-007&2.0246&2.0515e-006&2.0259\\
1/640&1.0397e-007&*&1.0397e-007&*&5.0375e-007&*\\
\hline
\end{tabular*}
\label{tab4}
\end{table}

\quad From Tables 3 and 4, we know that the numerical results are
consistent with our theoretical results.

\begin{table}[htbp]
\begin{tabular*}{14cm}{p{47pt}p{76pt}p{43pt}p{86pt}p{43pt}}
\multicolumn{5}{p{410pt}}{\small Table 5}\\
\multicolumn{5}{p{410pt}}{\small Errors of the numerical solutions in the mesh-grid points and all the points (new method) and the time to get those solutions when $\tau=h^2, t=1, h=1/(n-1)$ .}\\
\hline grid node number&\multicolumn{2}{p{76pt}}{H-OCD method}&\multicolumn{2}{p{86pt}}{numerical gradient scheme}\\
\cline{2-2}\cline{4-4}&Error&Time&Error&Time\\
\hline N=15&6.2975e-007&0.1248&1.3540e-005&0.0406\\
N=31&3.9376e-008&0.5725&1.1199e-006&0.1265\\
N=63&2.4630e-009&3.1590&8.0259e-008&0.5959\\
N=127&1.5453e-010&20.117&5.3654e-009&3.1844\\
N=255& 1.2050e-011& \fbox{128.44}&3.4669e-010&\fbox{21.542}\\
\hline
\end{tabular*}
\label{tab5}
\end{table}

\quad In addition, the conclusion in the space direction we get from Tab. 5 is the same as that
from Table 2. Thus, combining with Fig.4, the advantage of the numerical gradient scheme is obviously. In Table 6, we consider the Richardson extrapolation on this H-OCD scheme (\ref{eq7})-(\ref{eq9}) in time direction, the result is consistent with the Theorem 2.3.

\begin{table}[htbp]
\begin{tabular*}{13cm}{p{55pt}p{75pt}p{55pt}p{75pt}p{55pt}}
\multicolumn{5}{p{410pt}}{\small Table 6}\\
\multicolumn{5}{p{410pt}}{\small Errors and rate of all the points (new method) for Problem 4.1, 4.2 when $\tau=h, T=1$.}\\
\hline $\tau=h$ &\multicolumn{2}{p{60pt}}{Problem 4.1 } &\multicolumn{2}{p{60pt}}{Problem 4.2 }\\
\cline{2-2}\cline{4-4}&Error&$\frac{{error(h,\tau)}}{{error(h/2,\tau/2)}}$&Error&$\frac{{error(h,\tau)}}{{error(h/2,\tau/2)}}$\\
\hline
1/8&5.5147e-006&14.2567& 2.0369e-005&14.9661\\
1/16&3.8682e-007&15.6251& 1.3610e-006&15.4867\\
1/32&2.4756e-008&15.9066& 8.7881e-008&15.7438\\
1/64& 1.5563e-009&15.9741& 5.5819e-009&15.8711\\
1/128&9.7429e-011&15.9934&3.5170e-010&15.9209\\
1/256&6.0918e-012&*& 2.2091e-011&*\\
\hline
\end{tabular*}
\label{tab6}
\end{table}

\subsection{Numerical Experiments for the Two-dimensional Case}

\quad \textbf{Example 4.3} When
$$\begin{array}{l}
u\left( {x,y,0} \right) = \sin \left( {\pi x} \right)\sin \left( {\pi y} \right),\\
u\left( {0,y,t} \right) = u\left( {1,y,t} \right) = u\left( {x,0,t} \right) = u\left( {x,1,t} \right) = 0,
\end{array}$$
the exact solution of the problem (\ref{eq3.1}) is
\[
u\left( {x,y,t} \right) = {e^{ - 2{\pi ^2}t}}\sin \left( {\pi x} \right)\sin \left( {\pi y} \right),\left( {x,y,t} \right) \in \Omega  \times \left( {0,T} \right].
\]

\quad Next, let us observe and compare the numerical solutions from different methods.

\quad Table 7 lists the computational results of the mesh-grid points and intermediate points with different spacial step sizes
when time step size is fixed as $\tau=1/100000$. We can see that the convergence orders in space can reach $O(h^4)$ which is consistent with
the theoretical analysis (see Theorem 3.1-3.2) in this article. In addition, from Table 8, we see also that the numerical gradient scheme has the same the convergence order $O(\tau^2+h^4)$ as H-OCD method when the time and space step sizes are both changing.

\begin{table}[htbp]
\begin{tabular*}{14cm}{p{45pt}p{84pt}p{45pt}p{84pt}p{45pt}}
\multicolumn{5}{p{390pt}}{\small Table 7}\\
\multicolumn{5}{p{390pt}}{\small Errors and rate of intermediate points and numerical gradient $P$ (\ref{eq13}) in space direction with $\tau=1/100000$.}\\
\hline $h$ &\multicolumn{2}{p{129pt}}{H-OCD mesh-grid points}&\multicolumn{2}{p{129pt}}{intermediate points (New method)}\\
\cline{2-2}\cline{4-4}&Error&Rate&Error&Rate\\
\hline 1/4&5.3017e-011&3.9403&*&*\\
1/8&3.4536e-012&3.9878&4.6576e-012&3.9629\\
1/16&2.1769e-013&3.9805&2.9869e-013&3.9791\\
1/32&1.3791e-014&*&1.8939e-014&*\\
\hline
\end{tabular*}
\label{tab7}
\end{table}

\begin{table}[htbp]
\begin{tabular*}{13cm}{p{55pt}p{75pt}p{55pt}p{75pt}p{55pt}}
\multicolumn{5}{p{410pt}}{\small Table 8}\\
\multicolumn{5}{p{410pt}}{\small Errors and rate of all the points (new method) for Problem 4.3 when $\tau=h^2, h=1/{n-1},T=1$.}\\
\hline $N$ &\multicolumn{2}{p{60pt}}{H-OCD method} &\multicolumn{2}{p{60pt}}{numerical gradient}\\
\cline{2-2}\cline{4-4}&Error&$\frac{{error(h,\tau)}}{{error(h/2,\tau/2)}}$&Error&$\frac{{error(h,\tau)}}{{error(h/2,\tau/2)}}$\\
\hline
N=5&1.6485e-009&9.7908& 1.8257e-009&11.0802\\
N=10&1.6838e-010&15.6079& 1.6477e-010&15.3196\\
N=20&1.0788e-011&15.9751& 1.0755e-011&15.9015\\
N=40&6.7530e-013&*&6.7638e-013&*\\
\hline
\end{tabular*}
\label{tab8}
\end{table}

\begin{table}[htbp]
\begin{tabular*}{14cm}{p{45pt}p{68pt}p{55pt}p{45pt}p{75pt}p{55pt}}
\multicolumn{6}{p{410pt}}{\small Table 9}\\
\multicolumn{6}{p{410pt}}{\small A comparison of computation time between H-OCD method and numerical gradient scheme.}\\
\hline grid $\;$number &\multicolumn{2}{p{50pt}}{H-OCD method}&grid $\;$number&\multicolumn{2}{p{50pt}}{numerical gradient}\\
\cline{2-2}\cline{5-5}&Error&Time& &Error&Time\\
\hline n=16&1.6485e-009&0.0374&n=17&1.8257e-009&0.0421\\
n=81&1.6838e-010&0.4563&n=117&1.6838e-010&0.5756\\
n=224&3.3602e-011&2.8782&n=433&3.4081e-011&2.8860\\
n=361&1.0788e-011&9.0527&n=745&1.0788e-011&10.8556\\
n=624&4.4056e-012&25.2347&n=1233&4.4370e-012&26.8556\\
n=899&2.1338e-012&63.2506&\fbox{n=1783}&2.1347e-012&\fbox{66.8199}\\
\fbox{n=1599}&6.7530e-013&\fbox{422.8602}&n=3183&6.7638e-013&424.2073\\
\hline
\end{tabular*}
\label{tab9}
\end{table}

\begin{figure}\label{fig1}
  \centering
    \includegraphics[width=1.05\textwidth]{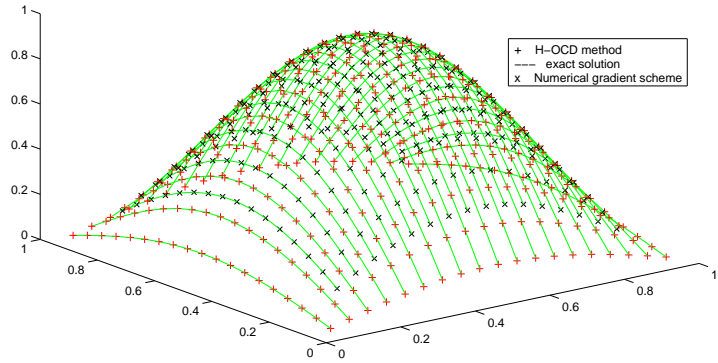}
  \caption{\small The curves of the numerical solution and exact solution for Example 4.3 when $h = 1/20,\tau  = 1/400, T= 1$.}
\end{figure}

In addition, Table 9 and Figure 6 also show similar results to those of Table 5 and Figure 3-4, respectively. Table 10 lists the computational results on Richardson extrapolation scheme. These results show that its convergence order in time direction can reach $O(\tau^4)$, which is consistent with
the theoretical analysis (see Theorem 3.3).

\begin{table}[htbp]
\begin{tabular*}{13cm}{p{55pt}p{75pt}p{55pt}p{75pt}p{55pt}}
\multicolumn{5}{p{410pt}}{\small Table 10}\\
\multicolumn{5}{p{410pt}}{\small The convergence order of Richardson extrapolation scheme for Problem 4.3 when $\tau=h/20, T=1$.}\\
\hline $h$ &\multicolumn{2}{p{60pt}}{H-OCD method} &\multicolumn{2}{p{60pt}}{Numerical gradient}\\
\cline{2-2}\cline{4-4}&Error&$\frac{{error(h,\tau)}}{{error(h/2,\tau/2)}}$&Error&$\frac{{error(h,\tau)}}{{error(h/2,\tau/2)}}$\\
\hline
h=1/5&2.1070e-011&14.1419&3.2895e-011&16.4508\\
h=1/10&1.4899e-012&15.9254& 1.9996e-012&15.7474\\
h=1/20&9.3555e-014&15.9822& 1.2698e-013&15.9392\\
h=1/40&5.8537e-015&15.9955&7.9665e-015&15.9848\\
h=1/40&3.6596e-016&*&4.9838e-016&*\\
\hline
\end{tabular*}
\label{tab10}
\end{table}

\quad For this two-dimension problem, we have obtained the similar experimental results as the previous one-dimension problem, which all shows that this method is effective.

\section{Conclusions}

\quad Recently, many people devote themselves on the development of
numerical approximation of heat equation problems. By the numerical
comparisons, we know that the high-order compact difference scheme (H-OCD) in \cite{Sun.Zhang} is better
than the traditional numerical schemes. In this article, we further improve this method to a new numerical
gradient scheme, which speeds up the convergence of the H-OCD scheme to some extent. Moreover, our theoretical analysis and numerical experiments show that this numerical gradient scheme has the same convergence order as H-OCD in \cite{Sun.Zhang}.

\textbf{Acknowledgements}. \emph{The authors sincerely thank the
reviewers and editor for their valuable and detailed comments and
suggestions on the early manuscript of this paper, which led to a
substantial improvement on the presentation and contents of this
paper.}

\end{document}